\documentclass[10pt,leqno]{article}
\usepackage{amsmath,amssymb}
\author{Rosa Cid-Mu$\tilde n$oz \and Manuel
Pedreira}
\title{Classification of Incidence Scrolls (II)}
\date{}

\newtheorem{teo}{Theorem}[section]

\newtheorem{prop}[teo]{Proposition}
\newtheorem{cor}[teo]{Corollary}
\newtheorem{lemma}[teo]{Lemma}

\newtheorem{rem}[teo]{Remark}
\def\di{\mathop{\rm dim}}
\newcommand{\p}{\ensuremath{\textrm{\rm I\hspace{-1.9pt}P}}}
\newcommand\om{\Omega}
\def\E{{\cal E}}
\def\Te{{\cal O}}

\font\euf=eufm10 at 14pt
\def\e{\mbox{\euf e}}
\def\d{\mbox{\euf b}}
\def\de{\mbox{\euf d}}
\def\b{{\cal B}}
\def\A{\alpha}

\def\qed{\hspace{\fill}$\rule{1.5mm}{1.5mm}$}
\def\lrw{\longrightarrow}
\def\rw{\rightarrow}
\def\sub{\subset}

\begin{document}
\maketitle

{\footnotesize {\bf Authors' address:} Departamento de
Algebra, Universidad de Santiago de Compostela. $15782$
Santiago de Compostela. Galicia. Spain. Phone:
34-81563100-ext.13152. Fax: 34-81597054. {\tt e}-mail: {\tt
rosacid@usc.es}; {\tt pedreira@zmat.usc.es}\\ {\bf
Abstract:} The aim of this paper is to obtain a classification of
the scrolls in $\p^n$ which are defined by a one-dimensional
family of lines meeting a certain set of linear spaces in
$\p^n$, a first classification for genus 0 and 1 is given in
paper \cite{rosa}. These ruled surfaces  are called {\it
incidence scrolls}, and such an indicated set is a {\it
base} of the incidence scroll. In particular, we compute its 
degree and genus. For this, we define the {\it fundamental
incidence scroll} to be the scroll in $\p^n$ formed by the
lines which meet $(2n-3)\, \p^{n-2}$'s in general position. 
Then all the others incidence scrolls will be portions of
degenerate forms of this. In this  way, we can obtain all
the incidence scrolls in $\p^n, \,n\geq3,$ with base in
general position.}
\vspace{0.1cm}

{\bf Introduction:} Throughout this paper, the base field for algebraic varieties
is ${\mathbb{C}}$. Let $\p^n$ be the n-dimensional complex
projective space and $G(l,n)$ the Grassmannian of $l$-planes
in
$\p^n$. Then $R^d_g \sub \p^n$ denotes a scroll of degree $d$
and genus $g$. We will follow the notation and terminology of
\cite{hatshor}.

It is useful to represent a scroll in $\p^n$ by a curve
$C\sub G(1,n)\sub \p^N$. The lines which intersect
a given subspace $\p^r\sub \p^n$ are represented by
the points of the special Schubert variety $\om (\p^r,\p^n)$.
Each
$\om (\p^r,\p^n)$ is the intersection of $G(1,n)$ with a
certain subspace of $\p^N$. Since $G(1,n)$ has dimension
$2n-2$ and we search a curve, we must impose $2n-3$ linear  
conditions on $G(1,n)$. Consequently, the choice of subspaces
is not arbitrary. Any set of subspaces of $\p^n$ which
imposes $2n-3$ linear conditions on $G(1,n)$ is the base of
a certain incidence scroll. The background about Schubert
varieties can be found in \cite{Kleiman}.

The aim of this paper is to obtain a classification of
the scrolls in
$\p^n$ which are defined by a one-dimensional family of
lines meeting a certain set of subspaces of $\p^n$, a first
classification for genus 0 and 1 is given in paper
\cite{rosa}. These ruled surfaces  are called {\it
incidence scrolls}, and such an indicated set is a {\it
base} of the incidence scroll. Unless otherwise stated, we
assume that the base spaces are in general position.

For the convenience of the reader we repeat the
relevant material from \cite{rosa} without proof, thus
making our exposition self-contained. Accordingly, Section 1
can be viewed as a brief summary of notation, fundamental
definitions and technical results which are used for our
analysis. Our first step will be to summarize some general 
properties of ruled surfaces and, in par\-ticular, of
incidence scrolls. Having revised the notion of incidence
scroll, we have compiled some basic properties of such a
scroll, the detailed proofs appear in \cite{rosa}. We
will expose a method to know when a scroll  is determined by
incidence. The affirmative solution would allow us to obtain
a base for each incidence scroll (there is really only one
way to choose this base). This is possible because the
families of directrix curves provide a natural and intrinsic
characterization of the incidence scrolls. The degree of the
scroll given by a general base is provided by Giambelli's
formula which appears in
\cite{G-h} or by Young tableaux (\cite{fulton}). Moreover,
the study of deformations of a given incidence scroll is a
powerful tool in order to simplify our proofs.
If the incidence scroll
$R^d_g
\sub \p^n$ breaks up into
$R^{d_{1}}_{g_{1}} \sub \p^r$ and $R^{d_{2}}_{g_{2}} \sub
\p^s$ with
$\delta$ generators in common, then $d=d_{1}+d_{2}$ and $
g=g_{1}+g_{2}+\delta-1$. 

In Section 2 we define the {\it fundamental incidence scroll}
to be the scroll in $\p^n$ formed by the lines which meet
$(2n-3)\,\p^{n-2}$ in general position. All the other
incidence scrolls will be portions of degenerate forms of
this. In this way, our main contribution is to establish
the degree and genus of the incidence scroll in $\p^n$
with base
$\b=\{i_s\,\p^{n-s-2},\ldots,i_1\,\p^{n-3},
(2n-3-\sum\limits_{j=1}^s(j+1)i_j\,\p^{n-2}\}$ in general
position, for any $i_1,\ldots,i_s$ between suitable limits.

Section 3 is devoted to study of decomposable ruled
surfaces. The important point to note here is the fact that
we find all the projective models of these which  are
defined by incidence. After having computed its degree and
genus, we restrict our analysis to those indecomposable
incidence scrolls in $\p^n$ whose two directrix curves of
minimum degree meet in a only point. Each of these is
projected from a point of a decomposable incidence scroll
keeping the incidence and genus.  

There are people who, even if not related to this paper,
have played a important role. My thanks to Proff. W.
Fulton, F. Sottile, I. Vainsencher and J. Kock. Special thanks
to Prof. A. Lascoux for useful advices (use Maple $6$ and
library ACE) and for precise bibliographic references. The
results on this paper  belong to the Ph.D. thesis of the
first author whose advisor  is the second one.

\section{Incidence Scrolls} \label{1}

A {\it ruled surface} is a surface $X$ together with a
surjective morphism $\pi \colon X\lrw C$ to a
smooth curve $C$ such that the fibre $X_y$
is isomorphic to $\p^1$ for every point $y\in C$, and such
that $\pi$ admits a section. There exists a locally free
sheaf $\E$ of rank 2 on $C$ such that 
$X\cong \p (\E)$ over $C$. Conversely, every such $\p
(\E)$ is a ruled surface over $C$. 

A {\it scroll} is a ruled surface embedded in $\p^n$ in such
a way that the fibres $f$ have degree $1$. If we take a very
ample divisor on $X$, $H \sim aC_o+\d f$, then the
embedding $\Phi_{|H|} \colon X
\lrw R^d_g\sub\p^n=\p(H^0(\Te_X(H)))$ determines a scroll
when $a=1$. A scroll $R^d_g\sub
\p^n$ is said to be {\it an incidence scroll} if it is
generated by the lines which meet a certain set $\b$ of
linear spaces in $\p^n$, or equivalently, if the
correspondent curve in
$G(1,n)$ is an intersection of special Schubert varieties
$\om(\p^{r}, \p^{n})$, $0\leq r <n-1$. Such a set is called
a {\it base} of the incidence scroll and such a base will be
denoted by:
$$\b =\{\p^{n_1},\p^{n_2},\cdots ,\p^{n_r}\}.$$ 
We will write it simply $\b$ when no confusion can arise,
where $n_1\leq n_2\leq \cdots \leq n_r$. 

Therefore,
unless otherwise stated, we will work with linear spaces in general position. By
general position we will mean that
$(\p^{n_1},\cdots ,\p^{n_r})\in {\cal W}=G(n_1,
n)\times \cdots \times G(n_r, n)$ is contained in
a nonempty open subset of
${\cal W}$. For simplicity of
notation, we abbreviate it to base in general
position. 
 
\begin{prop} \label{IS} ({\em \cite{rosa}, Proposition 2.4})
The intersection
$C=\bigcap_{i=1}^{r}
\Omega(\p^{n_i},\p^{n})$ of special Schubert varieties
associated to linear spaces $\p^{n_i},\, i=1,\cdots,r$, in
general position is an irreducible curve of
$G(1,n)$ if and only if it verifies the following
equality
$$ rn-(n_1+n_2+\cdots +n_r)-r=2n-3 \eqno(IS)$$\qed 
\end{prop} 

Moreover, the incidence scroll generated by a base $\b$
have degree $d$ if and only if we obtain the following equality of
Schubert cycles: $$\om(n_1, n) \cdots
\om(n_r, n)=d\, \om(0, 2).$$ Using Young tableaux, we can
also compute the degree of the scroll. Consider all ways of
filling a $2$ by $n-1$ rectangle with $n_1\, 1$'s, $n_2\,
2$'s, and so on until $n_r\, r$'s, and the one $r+1$ in the
lower right corner. The fillings have the restriction that
the numbers must be weakly increase across each row, but
strictly increase down each column. The number of such
objects is the degree that we desire (\cite{frank}). 

Furthermore, we present one of the three main theorems of
the symbolic formalism, known as Schubert calculus, for
solving enumerative problems.

\begin{teo} \label{pieri}
(Pieri's formula) For all sequences of integers $0\leq
a_0<\cdots <a_l \leq n$ and for $h=0,\cdots,n-l$, the
following formula holds in the cohomology ring
$H^{\star}(G(l,n); {\mathbb{Z}})$:
$$
\om(a_0,\cdots,a_l)\om(h, n)=\sum \om(b_0,\cdots,b_l)
$$
where the sum ranges over all sequences of integers
$b_0<\cdots <b_l$ satisfying $0\leq b_0\leq a_0<b_1\leq
a_1<\cdots<b_l\leq a_l$ and
$\sum^l_{i=0}b_i=\sum^l_{i=0}a_i-(n-l-h)$.
\end{teo}
{\bf Proof.} See \cite{Kleiman}, p. 1073.\qed

Finally, let us mention an important property of
degeneration of these scrolls.

\begin{prop} \label{degen} ({\em \cite{rosa}, Proposition
3.1}) Let $R^d_g \sub \p^n$ be an incidence scroll with base 
$\b$ in general position. Suppose that
$\p^{n_i}\vee\p^{n_j}=\p^{n-1}$ and have in common
$\p^{m}$, 
$m=n_i+n_j-n+1$. Then
the scroll breaks up into: 
\begin{enumerate} 
\item[-] $R^{d_1}_{g{_1}} \sub \p^{n}$
with base $\dot \b
=\{\p^{m},\p^{n_1},\cdots,\widehat{\p^{n_i}},\cdots,\widehat{\p^{n_j}},
\cdots,\p^{n_r}\}$ (which is possibly degenerate);
\item[-] $R^{d_2}_{g_2} \sub \p^{n-1}$
with base $\ddot \b =\{\cdots ,
\p^{n_{i-1}-1},\p^{n_i},\p^{n_{i+1}-1},\cdots
,\p^{n_{j-1}-1},$ $\p^{n_j},\p^{n_{j+1}-1},\cdots \}$  
\end{enumerate}
which have $\kappa \geq 1$ generators in common. Then,
$d=d_1+d_2$ and $g=g_1+g_2+ \kappa -1$.

Moreover, if $m=0$, then the incidence scroll breaks up
into a plane and an incidence scroll $R^{d-1}_g \sub
\p^{n-1}$ with base $\ddot \b$ in general position.\qed
\end{prop}

If $m=0$, then shall refer to this particular degeneration as
{\it join
$\p^{n_i}$ and $\p^{n_j}$} (i.e., $n_i+n_j=n-1$) and to the
inverse as {\it separate
$\p^{n_i}$ and $\p^{n_j}$} (i.e., $n_i+n_j=n$).

From now on, we will talk about a decomposable
incidence scroll if the corresponding ruled surface 
$X=\p(\E)$ has $\E$ decomposable. If $\E$ is indecomposable, 
then we will talk of an indecomposable incidence scroll.
 
Let $X=\p(\E)$ be a ruled surface
over the curve $C$ of genus $g$, determined by a
decomposable normalized bundle $\E\cong
\Te_C\oplus\Te_C(\e)$ such that $deg(-\e)=e\geq
0$. Let $H\sim C_o+\d f$ be the very ample divisor
on $X$ with $m=deg(\d)$ which gives the immersion
of the ruled surface as the scroll $R^d_g \sub
\p^n$ such that $d=2m-e$ and $n=2(m-g)-e+1+i$, being $i$ the 
speciality of the
scroll ($i=h^1(\Te_C(\d))+h^1(\Te_C(\d +\e))$).
Geometrically, $X$ has two disjoint directrix,
denoted by $C_o$ and $C_1$, such that $C_1\sim
C_o-\e f$. Moreover, it is easy to check that
these satisfy $\phi_{\d + \e }
\colon C_o\lrw C_g^{m-e} \sub \p^{m-e-g+i_1}$ and 
$\phi_{\d } \colon C_1\lrw C_g^{m}\sub
\p^{m-g+i_2}$, being $i_1=h^1(\Te_C(\d +\e))$ and
$i_2=h^1(\Te_C(\d))$. Therefore
$i=i_1+i_2$ and $\p^{m-e-g+i_1}\cap \p^{m-g+i_2} =
\emptyset$.

\begin{prop} \label{25} ({\em \cite{rosa}, Proposition 2.5})
If $R^d_g \sub \p^n$ is a decomposable incidence scroll with
base in general position, then
$\p^{m-e-g+i_1}$ and $\p^{m-g+i_2}$ are base spaces.\qed
\end{prop}

\begin{prop} \label{l62} ({\em \cite{rosa}, Proposition
2.7}) Let $R^{2m-e}_g\sub
\p^{2(m-g)-e+1+i}$ be a decomposable incidence
scroll with base in general position. 
\begin{enumerate}
\item[(a)] If $\e \not\sim 0$, then there are 
$(e+2-g+h^1(\Te_C(-\e)))\, \p^{m-g+i_2}$'s in
$\b$, whenever possible, i.e., when
$m-g+i_2+(e+2-g+h^1(\Te_C(-\e)))(m-e-g+i_1)\leq
4(m-g)-2e+2i-1$ and the generic curve in $|C_o-\e f|$ is
irreducible.
\item[(b)] If $\e \sim 0$, then there are 
$3\, \p^{m-e-g+i_1}$'s in
$\b$.\qed
\end{enumerate}
\end{prop}

A first classification of scrolls of genus $0$ and $1$ is
given in paper \cite{rosa}, Theorems 4.2 and 5.1. 

\begin{teo}\label{rael}
Let $X$ be a ruled surface over the curve
$C$ of genus $g=0,1$. Let
$H\sim C_o+\d f$ be a very ample divisor on
$X$ with $m=deg(\d)$ and let {\small $\Phi_H:X \rw
R_g^{2m-e}
\sub
\p^{2m-e-2g+1}$} be the closed immersion defined by
$|H|$. Then $R_g^{2m-e}$ is an incidence scroll
if and only if it satisfies one of the following
conditions:
\begin{enumerate}
\item $g=0$ and $e=0,1$;
\item $g=0$ and $m=e+1$;
\item $g=1, \, \,e=-1$ and $m=2$;
\item $g=1, \, \,\e\sim 0$ and $m=4$;
\item $X$ decomposable, $\e\not\sim 0,\,\, g=1, \, \,0 \leq
e\leq 3$ and $m=e+3$.\qed 
\end{enumerate}
\end{teo}

\begin{prop}\label{condicion} Let $X=\p(\E)$ be a
decomposable ruled surface over the curve $C$ of
genus $g$. Let
$H\sim C_0+\d f$ be a very ample divisor on $X$ which gives
the immersion of $X$ as the scroll
$R^{2m-e}_{g}\sub\p^{2(m-g)-e+i+1}$ such that $m=deg(\d)$. 
\begin{enumerate}
\item[1.]$(e-g+h^1(\Te_C(-\e)))(m-e-g+i_1)=m-g+i_2-1\,
(\star)$ and
$h^0(\Te_X(C_o-\e f))\leq 3 \Rightarrow R$ is an incidence
scroll with 
$\b=\{\p^{m-e-g+i_1}, 
(e+2-g+h^1(\Te_C(-\e))
)$ $\p^{m-g+i_2}\}$.
\item[2.] $(e-g+h^1(\Te_C(-\e)))(m-e-g+i_1)>
m-g+i_2-1\Rightarrow R$ is not
an incidence scroll. 
\end{enumerate}
\end{prop}
{\bf Proof}. $(\star)$ is
equivalent to
$(e-1-g+h^1(\Te_C(-\e)))(m-e-g+i_1)=e+i_2-i_1-1$. In
particular, $1\leq h^0(\Te_X(C_o-\e
f))\leq e+i_2-i_1+2$. 
\begin{enumerate}
\item[a.] For $h^0(\Te_X(C_o-\e f))=1$, we have
$2(m-g)-e+i=1$, which contradicts the
fact that $2(m-g)-e+i\geq2$.
\item[b.] For $h^0(\Te_X(C_o-\e f))=2$, we have
$m-e-g+i_1=m-g+i_2=1$. Then the incidence scroll is the
smooth quadric surface in $\p^3$ with base
$\b=\{3\,\p^1\}$.
\item[c.] For $h^0(\Te_X(C_o-\e f))=3$, i.e.,
$i_1-e=i_2-1$, we have $\b=\{\p^{m-g+i_2-1},$ $3\,
\p^{m-g+i_2}\}$. According to Theorem \ref{rael}, we know
that $\b$ generates a rational incidence scroll. The proof
is completed by showing that if $(\star)$ is true for
$g\geq 1$, then
$h^0(\Te_X(C_o-\e f))\not=3$. To do this, suppose that $g\geq
1$ and $h^0(\Te_X(C_o-\e f))=3$. Under the above conditions,
we see that
$e-1=g-h^1(\Te_C(-\e))$. Then can assume that  
$e\geq 1$, on the contrary ($e=0$), we find that
$h^1(\Te_C(-\e))$ is equal to
$g$ or $g-1$, hence that 
$e-1\not=g-h^1(\Te_C(-\e))$. Since $e\geq
1$, we take a general point $P$ of $R$ such that
$P\notin\p^{m-g+i_2-1}$ and $P\notin\p^{m-g+i_2}$ for some 
$\p^{m-g+i_2}\in\b$. The projection of 
$R^{2m-e}_{g}$ from $P$ is a scroll with at least
two directrix curves contain in linear spaces of dimension
$m-g+i_2-1$. Since $h^0(\Te_X(C_o-\e f))=3$, there is a
one-dimensional family of directrix curves $D\sim C_o-\e f$
that contain $P$. These curves
have the same degree, i.e., $m-e=m-1$. More precisely, this
means that $g=0$ because $R$ is a smooth irreducible scroll
of degree $2m-1$ in $\p^{2m}$ (see \cite{hatshor}; Remark
2.19.2).

For $g=0$, we know that $i_1=i_2=0$, hence we conclude from
$(\star)$ that $e=1$, and finally that $R$ is
rational of degree $2m-1$. It is known that 
$R^{2m-1}_{0}\sub
\p^{2m}$ is incidence with base
$\b=\{\p^{m-1},  3\,\p^{m}\}$ (\cite{rosa}).
\end{enumerate}

Next, suppose that 
$(e-g+h^1(\Te_C(-\e)))(m-e-g+i_1)>m-g+i_2-1$. For a suitable
number $\eta$ such that $1< \eta<m_1$, we have:
$$\begin{array}{l}
m-g+i_2+\eta(m-e-g+i_1)\leq4(m-g)-2e+2i-1\quad\quad\quad
(\star\star)\\ m-g+i_2+(\eta+1)(m-e-g+i_1)>4(m-g)-2e+2i-1.
\end{array}$$
Analysis similar to that in Proposition \ref{l62} shows that
if $R$ is defined by incidence, then there are
$\eta\,\p^{m-g+i_2}$'s in the base. Then $\b
=\{\p^{m-g-e+i_1},\eta\p^{m-g+i_2}, \p^{n_1},$ $\dots
,\p^{n_r}\}$ in general position with $m-g+i_2+1\leq n_i \leq
2(m-g)-e+i$ is the base of $R$. If $(\star\star)$ is
an inequality, then $n_1=\dots=n_r=2(m-g)-e+i$. Since the
scroll has at least one directrix curve $C^m_g\sub
\p^{m-g+i_2}\notin \b$ which is linearly independent from
the other directrix curves of degree $m$, we can take a
generic hyperplane $\p^{2(m-g)-e+i}$ through $\p^{m-g+i_2}$.
Then there are $(m-e)$ lines in $\p^{2(m-g)-e+i-\eta}$ which
meet 
$\p^{{n_1}-\eta-1},\dots,\p^{{n_r}-\eta-1}, \p^{m-e-g+i_1-1}$
and $\eta\,\p^{m-g+i_2-\eta}$, which is impossible.\qed

\begin{rem}{\em $(1.)$ is still true if we replace
$h^0(\Te_X(C_o-\e f))\leq 3$ by $g\leq 2$. For
$g=0$ (i.e., $i_1=i_2=0$) and $e\geq2$, we obtain
$m-e=1$ and
$h^0(\Te_X(C_o-\e f))=e+2$. Under these
conditions, we find $R^{e+2}_0\sub
\p^{e+3}$ with $\b=\{\p^1, (e+2)\, \p^{e+1}\}$ (see
Theorem \ref{rael}). The case $e=1$ has been study in the
above proof. For $g=1$, $(\star)$ is true only
for $e=3$ and $n=6$ (see Theorem \ref{rael}). Finally, for
$g(C)=2$, is trivial because $(\star)$ is impossible for
any $e\geq0$.}
\end{rem}

Finally, we will study the inequality
$(e-g+h^1(\Te_C(-\e)))(m-e-g+i_1)< m-g+i_2-1$. If $R$
is an incidence scroll, then there are  an $\p^{m-e-g+i_1}$
and $(e-g+2+h^1(\Te_C(-\e)))$
$\p^{m-g+i_2}$'s in $\b$. But these spaces are not
sufficient, so we must study the following directrix curves
of the scroll. These are curves of type $C^s_g\sub
\p^{s-g+i_s}$ such that $m-g+i_2+1\leq s-g+i_s\leq
2(m-g)-e+i$ which belong to a linear system $|C_o+(\de -\e
)f|$ where
$deg(\de)=s-m$. We proceed in three steps. 

\begin{enumerate}
\item[3.a)] If $s-g+i_s<2(m-g)-e+i-1$, then let $m_2\geq1$
be the number of directrix curves which are linearly
independent on $|C_o+(\de -\e )f|$.
\begin{enumerate}
\item[3.a.i)]
$(e-g+h^1(\Te_C(-\e)))(m-e-g+i_1)+m_2(2(m-g)-e+i-(s-g+i_s))<
m-g+i_2-1 \Rightarrow \{ R$ is an incidence scroll
$\Rightarrow
\{\p^{m-e-g+i_1},(e+2-g+h^1(\Te_C(-\e)))\,
\p^{m-g+i_2}, m_2\,\p^{s-g+i_s}\}\sub\b\} \Rightarrow$ to be
continued.
\item[3.a.ii)]$(e-g+h^1(\Te_C(-\e)))(m-e-g+i_1)+m_2(2(m-g)-e+i-(s-g+i_s))=
m-g+i_2-1 \Rightarrow \{ R$ is an incidence scroll 
$\Rightarrow
\b=\{\p^{m-e-g+i_1},(e+2-g+h^1(\Te_C(-\e)))\,
\p^{m-g+i_2},m_2\, \p^{s-g+i_s}\}\}$.
\item[3.a.iii)]
$(e-g+h^1(\Te_C(-\e)))(m-e-g+i_1)+m_2(2(m-g)-e+i-(s-g+i_s))>
m-g+i_2-1 \Rightarrow  \{ R$ is an incidence scroll
$\Leftrightarrow 
(e-g+h^1(\Te_C(-\e)))(m-e-g+i_1)+2(m-g)-e+i-(s-g+i_s)=
m-g+i_2-1 \}$ and  $\b=\{\p^{m-e-g+i_1},m_1\,
\p^{m-g+i_2}, \p^{s-g+i_s}\}$.
\end{enumerate}
\item[3.b)] $s-g+i_s=2(m-g)-e+i-1 \Rightarrow \{ R$ is
an incidence scroll $\Rightarrow
\b=\{\p^{m-e-g+i_1},(e+2-g+h^1(\Te_C(-\e)))\,
\p^{m-g+i_2},\eta \, \p^{s-g+i_s}\}$ with $\eta=m-g+i_2-1-
(e-g+h^1(\Te_C(-\e)))(m-e-g+i_1)\}$.
\item[3.c)] $s-g+i_s=2(m-g)-e+i \Rightarrow  R$ is not
an incidence scroll.
\end{enumerate}

\section{Fundamental Incidence Scroll}
\label{2}

We define {\it the fundamental incidence scroll} of
$\p^n$ to be the incidence scroll with base  $\b_n
=\{(2n-3)\,\p^{n-2}\}$, i.e., the curve
$C=\bigcap_{i=1}^{2n-3}\om_i(\p^{n-2}, \p^n)\sub G(1, n)$.
Using \cite {h-p}, p. 364, we deduce that it is a scroll of
degree
$\frac{1}{n-1}\binom{2n-2}{n}$ and all the other
incidence scrolls are portions of degenerate forms of this.
Moreover, we find that the degree of a
directrix curve contain in an $\p^{n-2}$ is
$3\frac{(2n-4)!}{(n-3)!n!}$. We have the quadric surface in
$\p^3$, elliptic quintic scroll generates by the lines which
meet $5\, \p^2$ in
$\p^4$ and so over.

Let $G:=G(1,n)$ be the Grassmannian of lines in a fixed
$\p^n$. The fixed $\p^n$ induces a trivial bundle 
$\Te_{G}^{n+1}$ on $G$ of rank $n+1$ and the $\p^1$
subspaces induce an subbundle $S$ of
$\Te_{G}^{n+1}$ of rank $2$. So on $G$
we have the universal exact sequence
$$0\lrw S\lrw \Te_{G}^{n+1}\lrw Q\lrw 0$$
where $Q$ is a vector bundle on $G$ of rank
$n-1$, called {\it the universal quotient bundle} on
$G$ and $S$ is called {\it the universal subbundle} on
$G$. This sequence determines a
homomorphism $S\lrw \Omega^1_G\bigotimes Q$. By duality,
it gives a homomorphism $T_G\lrw S^{\vee}\bigotimes Q$
which is an isomorphism. It follows easily that $$0\lrw
S\bigotimes S^{\vee}\lrw
\Te_{G}^{n+1}\bigotimes S^{\vee}\lrw
Q\bigotimes S^{\vee}\lrw 0,$$ hence that
$\wedge^{n+1}
T_G=\wedge^{n+1}(\Te_{G}(1)^{n+1})=
\bigotimes_{1}^{n+1}\Te_{G}(1)=\Te_{G}(n+1)$,
and finally that $\Te_G(K_G)=\Te_G(-n-1)$, where $K_G$
denotes the canonical divisor on $G$.
Let us apply this to
$C$ which is the intersection of $G$ with $2n-3$ generic
hyperplanes of $\p^N$. By
\cite{hatshor}; II, Exercise 8.4, we obtain $\Te_C(K_C)\cong
\Te_C(n-4)$ and hence
$$g_C=\frac{n-4}{2n-2}\binom{2n-2}{n}-1.$$

\begin{prop}\label{fun} The fundamental incidence scroll of
$\p^n$ has degree
$d(n)=\frac{1}{2}\binom{2n-2}{n}$ and
genus $g(n)=\frac{n-4}{2n-2}\binom{2n-2}{n}-1$. Moreover, the
degree of the minimum directrix curve is
$3\frac{(2n-4)!}{(n-3)!n!}$.\qed
\end{prop}

If $n\geq 5$, then we obtain an indecomposable incidence
scroll because $\di$ $(\p^{n-2}\cap\p^{n-2})=n-4$. If we
take two minimum directrix curves $C\sub\p^{n-2}$ and
$\tilde{C}\sub
\p^{n-2}$, then $C \cdot \tilde{C}
=2(n-3)\frac{(2n-4)!}{(n-2)!n!}=-e$.

\begin{prop}\label{noun} Let
$R^{d(n,i)}_{g(n,i)}\sub
\p^n$ be an incidence scroll with base
$\b(n,i)=\{i\,\p^{n-3},(2n-3-2i)\,\p^{n-2}\}$ in
general position. Then:
\begin{enumerate}
\item
$d(n,i)=\sum\limits_{k=0}^{i}(-1)^k\binom{i}{k}d(n-k)$,
\item
$g(n,i)=\sum\limits_{k=0}^{i}(-1)^k\binom{i}{k}g(n-k)-
\sum\limits_{k_1=0}^{i-1}\sum\limits_{k_2=0}^{k_1}(-1)^{k_2}
\binom{k_1}{k_2}d'(n-k_2-1, i-k_1-1)$.
\end{enumerate}
\end{prop}
{\bf Proof.} The proof is by induction on $i\geq 0$.
By the above proposition, $d(n,0)=d(n)$ and $g(n,0)=g(n)$.
Suppose that the incidence scroll with
base $\b(n,i-1)$
in general position has degree $d(n,i-1)$ and genus
$g(n,i-1)$. Since $2n-3-2(i-1)>2$, we can suppose that
$2\,\p^{n-2}$'s are contained in a hyperplane,
then $\b(n,i-1)$ breaks up into
$\b(n,i)\sub\p^n$ and $\b(n-1,i-1)\sub\p^{n-1}$ with
$d(n-1,i-1)$ generators in common. Whence,
\begin{enumerate}
\item $d(n,i)=d(n,i-1)-d(n-1,i-1)=
\sum\limits_{k=0}^{i-1}
(-1)^{k}\binom{i-1}{k}d(n-k)-\sum\limits_{k=0}^{i-1}
(-1)^{k}\\
\binom{i-1}{k}d(n-k-1)=
\sum\limits_{k_1=0}^{i}
(-1)^{k}\binom{i}{k}d(n-k)$,
\item
$g(n,i)=g(n,i-1)-g(n-1,i-1)-d(n-1,i-1)+1=\sum\limits_{k=0}^{i}
(-1)^{k}\binom{i}{k}g(n-k)-\sum\limits_{k_1=0}^{i-1}\sum\limits_{k_2=0}^{k_1}(-1)^{k_2}
\binom{k_1}{k_2}[d(n-k_2-1, i-k_1-1)-1].$\qed
\end{enumerate}

\begin{teo}\label{todas} Let
$R^{d(n,i_1,\ldots,i_s)}_{g(n,i_1,\ldots,i_s)}\sub
\p^n$ be an incidence scroll with base
$\b(n,i_1,$ $\ldots,i_s)$
$=\{i_s\,\p^{n-s-2},\ldots,i_1\,\p^{n-3},
(2n-3-\sum_{j=1}^{s}(j+1)i_j)\,\p^{n-2}\}$ in general
position. Then
\begin{enumerate}
\item{\small $d(n,i_1,\ldots,i_s)=\sum\limits_{k_1=0}^{i_s}
\sum\limits_{k_2=0}^{i_s+i_{(s-1)}-k_1}
\ldots
\sum\limits_{k_s=0}
^{\sum\limits_{j=1}^{s}i_j-k_{(s-1)}}(-1)^
{\sum\limits_{j=1}^{s}k_j}\binom{i_s}{k_1}
\binom{i_s+i_{(s-1)}-k_1}{k_2}\\
\ldots
\binom{\sum\limits_{j=1}^{s}i_j-k_{(s-1)}}{k_{s}}
d(n-\sum\limits_{j=1}^{s}k_j),$} 
\item{\small $g(n,i_1,\ldots,i_s)=\sum\limits_{k_1=0}^{i_s}
\sum\limits_{k_2=0}^{i_s+i_{(s-1)}-k_1}
\ldots
\sum\limits_{k_s=0}
^{\sum\limits_{j=1}^{s}i_j-k_{(s-1)}}(-1)^
{\sum\limits_{j=1}^{s}k_j}\binom{i_s}{k_1}
\binom{i_s+i_{(s-1)}-k_1}{k_2}\\
\ldots
\binom{\sum\limits_{j=1}^{s}i_j-k_{(s-1)}}{k_{s}}g(n-\sum\limits_{j=1}^{s}k_j)
-\sum\limits_{k_1=0}^{i_s-1}\sum\limits_{k_2=0}^{k_1}
(-1)^{k_2}
\binom{k_1}{k_2} \bigl[ d(n-k_2-1,i_1,\ldots,i_{(s-3)},\\
i_{(s-2)}
+k_2,i_{(s-1)}+1+k_1-k_2,
i_s-k_1-1)-1 \bigr]-\sum\limits_{\A=3}^{s+1}\Bigl[
\sum\limits_{k_1=0}^{i_s}
\sum\limits_{k_2=0}^{i_s+i_{(s-1)}-k_1}
\ldots\\
\sum\limits_{k_{(\A-2)}=0}
^{\sum\limits_{j=s-\A+3}^{s}i_j-k_{(\A-3)}}
\sum\limits_{k_{(\A-1)}=0}
^{\sum\limits_{j=s-\A+2}^{s}i_j-k_{(\A-2)}-1}
\sum\limits_{k_{\A}=0}
^{k_{(\A-1)}}
(-1)^{\sum\limits_{j\not=\A-1}k_j}\binom{i_s}{k_1}
\binom{i_s+i_{(s-1)}-k_1}{k_2}
\ldots\\
\binom{\sum\limits_{j=s-\A+3}^{s}i_j-k_{(\A-3)}}{k_{(\A-2)}}
\binom{k_{(\A-1)}}{k_{\A}} \bigl[ d(n-1
-\sum\limits_{j\not=\A-1}k_j,
i_1,\ldots,i_{(s-\A-1)},
i_{(s-\A)}+
+k_{(\A-1)},\\i_{(s-\A+1)}+1+k_{(\A-2)}+k_{(\A-1)}
-k_{\A},\sum\limits_{j=s-\A+2}^{s}i_j-k_{(\A-2)}-k_{(\A-1)}-1)-1\bigr]
\Bigr]$.}
\end{enumerate}

\end{teo}
{\bf Proof.} We proceed by induction on
$s\geq 1$. The case $s=1$ is Proposition \ref{noun}. If
$s=2$, we will prove that the incidence scroll with base
$\b(n,i_1,i_2)\sub
\p^n$ has 
\begin{enumerate}
\item
$d(n,i_1,i_2)=\sum\limits_{k_1=0}^{i_2}\sum\limits_{k_2=0}
^{i_1+i_2-k_1}(-1)^{k_1+k_2}\binom{i_2}{k_1}
\binom{i_1+i_2-k_1}{k_2}d(n-k_1-k_2)$,
\item
$g(n,i_1,i_2)=\sum\limits_{k_1=0}^{i_2}\sum\limits_{k_2=0}
^{i_1+i_2-k_1}(-1)^{k_1+k_2}\binom{i_2}{k_1}
\binom{i_1+i_2-k_1}{k_2}g(n-k_1-k_2)-
\sum\limits_{k_1=0}^{i_2-1}\sum\limits_{k_2=0}^{k_1}\\
(-1)^{k_2}
\binom{k_1}{k_2}\bigl[d(n-k_2-1,i_1+1+k_1-k_2,
i_2-k_1-1)-1\bigr]-
\sum\limits_{k_1=0}^{i_2}\sum\limits_{k_2=0}^
{i_1+i_2-k_1-1}\sum\limits_{k_3=0}^{k_2}\\
(-1)^{k_1+k_3}
\binom{i_2}{k_1}\binom{k_2}{k_3}[d(n-k_3-k_1-1,i_1+i_2-1-k_1-k_2)-1]$.
\end{enumerate} 
To do this, we apply induction on $i_2\geq 0$. It 
is true for $i_2=0$, by Proposition \ref{noun}. Supposing
the theorem true for $i_2-1\geq0$, we prove it for any $i_2$.
Since $(2n-3-2i_1-3i_2)\geq 0$, 
$(2n-3-2i_1-3(i_2-1))\geq 3$. Let us suppose for the moment
$i_1\geq 1$. Then we can consider that
$\p^{n-3}\vee\p^{n-2}=\p^{n-1}$, for any
$\p^{n-3}, \p^{n-2} \in \b(n,i_1,i_2-1)$. Then
$\b(n,i_1,i_2-1)$ degenerates into
$\b(n,i_1-1,i_2)\sub\p^n$ and
$\b(n-1,i_1-1,i_2-1)\sub\p^{n-1}$ with 
$d(n-1,i_1,i_2-1)$ generators in common. Whence,
writing $i_1+1$ instead $i_1$, we find that
\begin{enumerate}
\item $d(n,i_1,i_2)=d(n,i_1+1,i_2-1)-d(n-1,i_1,i_2-1)
=\sum\limits_{k=0}^{i_2}
(-1)^{k}\binom{i_2}{k}d(n-k,i_1+i_2-k)$,
\item
$g(n,i_1,i_2)=g(n,i_1+1,i_2-1)-g(n-1,i_1,i_2-1)-
d(n-1,i_1+1,i_2-1)+1=\sum\limits_{k=0}^{i_2}
(-1)^{k}\binom{i_2}{k}g(n-k,i_1+i_2-k)-
\sum\limits_{k_1=0}^{i_2-1}
\sum\limits_{k_2=0}^{k_1}(-1)^{k_2}
\binom{k_1}{k_2}[d(n-k_2-1,i_1+1+k_1-k_2,i_2-k_1-1)-1].$
\end{enumerate}
For $i_1=0$, we can suppose that 
$2\,\p^{n-2}$'s of $\b(n,0,i_2-1)$ are containing in a
hyperplane. Then we obtain $\b(n,1,i_2-1)$ which verifies
the above assumptions.

Assume the theorem holds for $s-1\geq 2$. For $s$, we use
induction on $i_s\geq 0$. An argument similar to $s=2$ shows
that the incidence scroll with base $\b(n,i_1,\ldots,i_s)$
has 
\begin{enumerate}
\item
$d(n,i_1,\ldots,i_s)=d(n,i_1,\ldots,i_{(s-1)}+1,i_s-1)
-d(n,i_1,\ldots,i_{(s-2)}+1,\\i_{(s-1)},i_s-1)=\sum\limits_{k=0}^{i_s}
(-1)^{k}\binom{i_s}{k}d(n-k,i_1,\ldots,i_{(s-3)},
i_{(s-2)}+k,i_{(s-1)}+i_s-k)$,
\item $g(n,i_1,\ldots,i_s)=g(n,i_1,\ldots,i_{(s-1)}+1,i_s-1)
-g(n,i_1,\ldots,i_{(s-2)}+1,\\i_{(s-1)},i_s-1)-
d(n,i_1,\ldots,i_{(s-1)}+1,i_s-1)+1=\sum\limits_{k=0}^{i_s}
(-1)^{k}\binom{i_s}{k}g(n-k,i_1,\ldots,i_{(s-3)},
i_{(s-2)}+k,i_{(s-1)}+i_s-k)-
\sum\limits_{k_1=0}^{i_s-1}
\sum\limits_{k_2=0}^{k_1}(-1)^{k_2}
\binom{k_1}{k_2}[d(n-k_2-1,i_1,\ldots,i_{(s-3)},
i_{(s-2)}+k_2,i_{(s-1)}+1+k_1-k_2,i_s-k_1-1)-1].$\qed
\end{enumerate}
\begin{rem}{\em When we talk about $\b(n,i_1,\ldots,i_s)$,
we must assume that $i_j\geq 0$ for $1\leq j\leq s$ and
$2n-3-\sum_{j=1}^{4}(j+1)i_j\geq 0$. Let us mention two
important consequences and, in particular,
the last extends the results of \cite{rosa2}. }
\end{rem}

\begin{cor}\label{nodescomponibles} Let $\b$ be defined by
$\b(n,i_1,\ldots,i_{(r-2)},0,\ldots,0,i_{(n-r-2)}=1)$ in
general position such that $2\leq r\leq n-2$. Then
$\b$ is the base of an incidence scroll of
degree $d(n,i_1,\ldots,i_{(r-2)},0,\ldots,0,i_{(n-r-2)}=1)$
and genus
$g(n,i_1,\ldots,i_{(r-2)},0,\ldots,0,i_{(n-r-2)}=1)$. Moreover,
these are all indecomposable incidence scrolls of
$\p^n$.\qed
\end{cor}

\begin{cor}
$\b^n_{i_1,\ldots,i_{(r-1)}}:=\b(n.i_1,\ldots,i_{(r-1)},0,
\ldots,0,i_{(n-r-2)}=1)$ in general position generates the
incidence scroll of degree
$d^n_{i_1,\ldots,i_{(r-1)}}=d(n,i_1,\ldots,$ $i_{(r-1)},0,$
$\ldots,0,$ $i_{(n-r-2)}=1)$ and genus
$g^n_{i_1,\ldots,i_{(r-1)}}=g(n,i_1,\ldots,i_{(r-1)},0,\ldots,$ 
$0, i_{(n-r-2)}$ $=1)$. Furthermore, these are all
incidence scrolls with base $\p^r$.\qed
\end{cor}

\section{Decomposable Incidence Scrolls}
\label{3}

Let $\pi \colon X=\p(\E)\lrw C$ be a 
geometrically ruled surface over a curve $C$ of genus
$g$ determined by a decomposable normalized bundle $\E \cong
\Te_C\oplus \Te_C(\e )$. Let $H\sim C_o+\d f$ be a very ample
divisor on $X$ with $m=deg(\d)$ which gives the immersion of
$X$ as the scroll $R^d_g\sub \p^n$. Following the above
notation, we will divided the study in two case.

\begin{rem}{\em In the sequel, when $\sum$ has not 
meaning, i.e., the lower limit is bigger than the higher
limit, we adopt the convention that it is $0$. Moreover, we
will write
$d(r,h_1,\ldots,h_s), g(r,h_1,\ldots,h_s),
\b(r,h_1,\ldots,h_s)$ and
$\Delta(r,h_1,\ldots,h_s)$, instead of
$d_{e}(r,h_1,\ldots,h_s), g_e(r,h_1,\ldots,h_s),
\b_e(r,h_1,\ldots,h_s)$ and
$\Delta _e(r,h_1,\ldots,$ $h_s)$ respectively, when no
confusion can arise. Finally, for abbreviation, we will write
$\Delta'(r,h_1,\ldots,h_s):=\Delta(r,h_1,\ldots,h_s)-1$.}
\end{rem}

\subsection{Incidence Scrolls with $\e\sim
0$}\label{inesp}

Suppose $\e\sim 0$. If $R^{2m}_g\sub \p^{2(m-g)+i+1}$ is an
incidence scroll, then there are $3\,\p^{m-g+i_1}$'s in $\b$
which impose $3(m-g+i_1)$ independent conditions
on $G(1,2(m-g+i_1)+1)$. Writing $r$ instead of $(m-g+i_1)$,
we may set up a one-to-one correspondence between incidence
scrolls $R\sub\p^{2r+1}$ with $\e\sim 0$ and partitions of
$r-1$, for any $r\geq 1$. The {\it partition} of a positive
integer $r$ is, by definition, a sequence of weakly
decreasing positive integers $\lambda =(\lambda _1 \geq
\ldots \geq\lambda_s)$ such that
$\sum\limits_{i=0}^s\lambda_i=r$.
 
\begin{prop} The number of incidence scroll
$R\sub\p^{2r+1}$ with
$\e\sim 0$ is exactly the number of partitions of $r-1$.
\end{prop}
{\bf Proof} Let $\lambda =(\lambda _1 \geq
\ldots \geq\lambda_s)$ be a partition of $r-1$. For each
$\lambda _i, \, 1\leq i\leq s$, we take a
generic $\p^{2r-\lambda _i}\sub\p^{2r+1}$ as base space. 
Hence every $\lambda$ gives rise an incidence scroll of
$\p^{2r+1}$ with base $\b=\{3\,\p^r, \p^{2r-\lambda _1},
\cdots,
\p^{2r-\lambda _s}\}$ where $\e\sim 0$. Conversely, to
every $\b_{\lambda}=\{3\,\p^r, \p^{n_1},
\ldots,\p^{n_s}\}\sub\p^{2r+1}$ such that $r+1\leq n_1$
there corresponds a partition $(2r-n_1\geq
\ldots\geq 2r-n_s)$ of $r-1$.\qed
 
\begin{lemma}\label{delta} For each partition $(r-h_1 \geq
h_1-h_2\geq \ldots \geq h_{(s-1)}-h_s \geq
1,{\stackrel{(h_s-1)}\ldots },1)$ of
$r-1$, we define $\Delta_{\e\sim
0}(r,h_1,\ldots,h_s):=\sigma(r)^3\cdot 
\sigma(r+h_1-1)\cdot\sigma(2r-h_1+h_2)\cdot\ldots\cdot\sigma(2r-h_{(s-1)}+h_s)
\cdot\sigma(2r-1)^{(h_s-1)}$ an intersection of Schubert
cycles in $\p^{2r+1}$. Then:
\begin{enumerate}
\item $\Delta_{\e\sim
0}(r,h)=\sum\limits_{k=0}^{h-1}\binom{h-1}{k}(r-2k+1)$;
\item $\Delta_{\e\sim 0}(r,h_1,\ldots,h_s)=
\sum\limits_{k_1=0}^{h_s-1}\sum\limits_{k_2=0}^
{h_{(s-1)}-h_s}
\sum\limits_{k_3=0}^
{h_{(s-2)}-h_{s}-2k_2}\ldots\sum\limits_{k_{(s-1)}=0}^
{h_2-h_s-2\sum\limits_{i=2}^{s-2}k_i}\\
\binom{h_s-1}{k_1}
(h_1-h_s+1-2\sum\limits_{j=2}^{s-1}k_j)(r-2k_1-h_1+h_s+1)$,
for  $s\geq 2.$
\end{enumerate}
\end{lemma}
{\bf Proof.} $1.$ For $s=1$, we apply induction on
$h\geq1$. Then $\Delta(r,1)=r+1$ because it is
the degree of a directrix curve, which is containing in an
$\p^{r+1}$, of the incidence scroll $R^{2r}_0\sub
\p^{2r+1}$. Suppose that the lemma is true for $h-1\geq1$,
we prove it for any $h$. From Pieri's formula, we have
{\small
$$\Delta(r,h)=\Delta(r,h-1)+\Delta(r-2,h-1)=
\sum\limits_{k=0}^{h-1}\binom{h-1}{k}\Delta(r-2k,1).$$} The
proof of the last equality is by induction on $h\geq1$,
using the fact that
$\binom{h-1}{k}=\binom{h-2}{k}+\binom{h-2}{k-1}$.

$2.$ Suppose $s=2$. We will see in Theorem \ref{emo} that
$\Delta(r,h_1,1)=h_1(r-h_1+2)$ because the
incidence scroll which corresponds to the partition
$(r-h_1,h_1-1)$ is a nonspecial scroll of degree
$2h_1(r-h_1+1)$ and genus $(r-h_1)(h_1-1)$. In general,
using Pieri's formula, {\small
$\Delta(r,h_1,h_2)=\Delta(r,h_1-1,h_2-1)+\Delta(r-2,h_1-1,h_2-1)=\sum\limits_{k=0}^{h_2-1}\binom{h_2-1}{k}
\Delta(r-2k,h_1-h_2+1,1)=\sum\limits_{k=0}^{h_2-1}\binom{h_2-1}{k}
(h_1-h_2+1)(r-2k-h_1+h_2+1).$}

Assume the formula holds for $s-1\geq 1$; we will prove it
for $s$. We have divided the proof into $2$ steps. 
\begin{enumerate}
\item[(a)] For $h_s=1$, according to Pieri's formula,
{\small
$\Delta(r,h_1,\ldots,h_{(s-1)},1)=\sum\limits_{k=0}^{h_{(s-1)}-1}\\
\sigma(r-2k)^3 
\sigma(r+h_1-4k-1)\sigma(2r-h_1+h_2-4k)\ldots
\sigma(2r-h_{(s-3)}+h_{(s-2)}-4k)
\sigma(2r-h_{(s-2)}-2k+1))_{|_{2r-4k+1}}$ $=\sum\limits_{k=0}^
{h_{(s-1)}-1}\Delta(r-2k,h_1-2k,\ldots,h_{(s-2)}-2k,1)\underset{\mbox{\tiny
induction on}\, s-1}{=}
\sum\limits_{k_1=0}^{h_{(s-1)}-1}
\sum\limits_{k_2=0}^{h_{(s-2)}-1-2k_1}
\ldots
\sum\limits_{k_{(s-2)}=0}^{h_2-1-2\sum\limits_{j=1}^{s-3}k_j}
(h_1-2\sum\limits_{j=1}^{s-2}k_j)(r-h_1+2)$}. The notation
$_{|_k}$ means that the intersection is in $\p^{k}$
instead of $\p^{2r+1}$.
\item[(b)] For $h_s\geq2$,
{\small
$\Delta(r,h_1,\ldots,h_{s})=\Delta(r,h_1-1,\ldots,h_{s}-1)+\Delta
(r-2,h_1-1,\ldots,h_{s}-1)$.} Then, for any
$h_s\geq1$, {\small
$\Delta(r,h_1,\ldots,h_{s})=\sum\limits_{k=0}^
{h_{s}-1}\binom{h_s-1}{k}\Delta(r-2k,h_1-h_s+1,
\ldots,h_{(s-1)}-h_s+1,1)$.}\qed
\end{enumerate}

\begin{prop} Each partition
$(r-h,1^{(h-1)}): =(r-h,1,{\stackrel{(h-1)} \ldots },1)$ of
$r-1$ gives rise to an incidence scroll $R^{d_{\e\sim
0}(r,h)}_{g_{\e\sim 0}(r,h)}\sub \p^{2r+1}$ with base
$\b_{\e\sim 0}(r,h)=\{3\,\p^r,  \p^{r+h},
(h-1)\,
\p^{2r-1}\}$ in general position, such that 
\begin{itemize}
\item $d_{\e
\sim 0}(r,h)=2\sum\limits_{k=0}^{h-1}\binom{h-1}{k}(r-2k),$
\item $g_{ \e \sim
0}(r,h)=\sum\limits_{k_1=2}^{h}\sum\limits_{k_2=0}^{h-k_1}
\binom{h-k_1}{k_2}
\Delta'_{ \e \sim
0}(r-2k_2-2,k_1-1).$
\end{itemize}

\end{prop}
{\bf Proof}. From Theorem \ref{rael}, it follows
that $d(r,1)=2r$ and $g(r,1)=0$. Suppose now that $h>1$. If,
in $\b(r,h)$, 
$\p^{r+h}\vee \p^{2r-1}=\p^{2r}$, then $R$ breaks up into
$R^{d(r-2,h-1)}_{g(r-2,h-1)}\sub
\p^{2r-3}$ with $\b(r-2,h-1)=\{3\,\p^{r-2},  \p^{r+h-3},
(h-2)\, \p^{2r-5}\}$ and $R^{d(r,h-1)}_{g(r,h-1)}\sub
\p^{2r+1}$ with $\b(r,h-1)=\{3\,\p^r,  \p^{r+h-1},
(h-2)\, \p^{2r-1}\}$, and $\Delta(r-2,h-1)$
generators in common. Whence, by induction on $h\geq1$,
\begin{enumerate}
\item$d(r,h)=d(r-2,h-1)+d(r,h-1)=\sum\limits_{k=0}^{h-1}\binom{h-1}{k}
d(r-2k,1),$
\item$g(r,h)=g(r-2,h-1)+g(r,h-1)+\Delta'(r-2,h-1)=
\sum\limits_{k_1=2}^{h}\sum\limits_{k_2=0}^{h-k_1}
\binom{h-k_1}{k_2}\Delta'(r-2k_2-2,k_1-1)$.\qed
\end{enumerate}

\begin{teo}\label{emo} For $s\geq2$, each partition
$(r-h_1,h_1-h_2,\ldots,h_{(s-1)}-h_s,1^{(h_s-1)})$ 
of $r-1$ gives rise to $R^{d_{\e\sim
0}(r,h_1,\ldots,h_s)}_{g_{\e\sim
0}(r,h_1,\ldots,h_s)}\sub\p^{2r+1}$ with 
$\b_{\e\sim 0}(r,h_1,
\ldots,h_s)=\{3\,\p^r, 
\p^{r+h_1},\p^{2r-h_1+h_2},\ldots, \p^{2r-h_{(s-1)}+h_s},
(h_s-1)
\p^{2r-1}\}$ in general position, such that
\begin{itemize}
\item {\small $d_{\e\sim 0}(r,h_1,\ldots,h_s)=2(h_1-h_2+1)
\sum\limits_{k_1=0}^{h_s-1}\sum\limits_{k_2=0}^
{h_{(s-1)}-h_s}\ldots\sum\limits_{k_{(s-1)}=0}^{h_2-h_3}
\binom{h_s-1}{k_1}
(r-h_1+h_2-2\sum\limits_{i=1}^{s-1}k_i)$,}
\item {\small $g_{\e\sim 0}(r,h_1,\ldots,h_s)=(h_1-h_2)
\sum\limits_{k_1=0}^{h_s-1}\sum\limits_{k_2=0}^
{h_{(s-1)}-h_s}\ldots\sum\limits_{k_{(s-1)}=0}^{h_2-h_3}
\binom{h_s-1}{k_1}
(r-h_1+h_2-1-2\sum\limits_{i=1}^{s-1}k_i)+
\sum\limits_{k_1=2}^{h_s}
\sum\limits_{k_2=0}^{h_s-k_1}\binom{h_s-k_1}{k_2}
\Delta'_{\e\sim0}(r-2k_2-2,
h_1-h_s+k_1-1,\ldots,h_{(s-1)}-h_s+k_1-1,k_1-1)+
\sum\limits_{\A=2}^{s-1}\Bigl[\;
\sum\limits_{k_1=0}^{h_s-1}
\sum\limits_{k_2=0}^{h_{(s-1)}-h_s}\ldots
\sum\limits_{k_{(s-\A)}=0}^{h_{(\A+1)}-h_{(\A+2)}}
\sum\limits_{k_{(s-\A+1)}=0}^{h_{\A}-h_{(\A+1)}-1}\binom{h_s-1}{k_1}
\Delta'_{\e\sim0}(r-2-2\sum\limits_{i=1}^{s-\A+1}k_i, 
h_1-h_{\A}+1,\ldots,
h_{(\A-1)}-h_{\A}+1,1)\Bigr].$}
\end{itemize}
\end{teo}
{\bf Proof.} Suppose $h_s=1$. We first compute $s=2$.
Therefore, $d(r,2,1)=4r-4$ and $g(r,2,1)=r-2$, by
\cite{rosa}, Example 3.3. By induction on
$h_1\geq 1$, we have
\begin{itemize} 
\item {\small
$d(r,h,1)=d(r,1)+d(r-2,h-1,1)=2r+2(h-1)(r-h)=2h(r-h+1);$}
\item {\small $g(r,h,1)
=g(r,1)+g(r-2,h-1,1)+r-2=(r-1-h)(h-2)+(r-2)=(r-h)(h-1).$}
\end{itemize}
Assuming theorem to hold for $s-1$, we will prove it
for $s$. To this end, we suppose that, in
$\b(r,h_1,\ldots,h_{(s-1)},1)$,
$\p^{r+h_1}\vee\p^{2r-h_{(s-1)}+1}=\p^{2r}$. Then we obtain
$\b(r,h_1-h_{(s-1)}+1,\ldots,h_{(s-2)}-h_{(s-1)}+1,1)\sub\p^{2r+1}$
and $\b(r-2,h_1-1,\ldots,h_{(s-1)}-1,1)\sub\p^{2r-3}$
with
$\Delta(r-2,h_1-h_{(s-1)}+1,\ldots,h_{(s-2)}-h_{(s-1)}+1,1)$
generators in common. In particular, using that
$\b(r,h_1,\ldots,h_{(s-2)},1,1)=\b(r,h_1,\ldots,h_{(s-2)},1)$,
we deduce that
\begin{itemize}
\item {\small $d(r,h_1,\ldots,h_{(s-1)},1)=
d(r,h_1-h_{(s-1)}+1,\ldots,h_{(s-2)}-h_{(s-1)}+1,1)
+d(r-2,h_1-1,\ldots,h_{(s-1)}-1,1)=\sum\limits_{k=0}^{h_{(s-1)}-1}
d(r-2k,h_1-h_{(s-1)}+1,\ldots,
h_{(s-2)}-h_{(s-1)}+1,1)=2(h_1-h_2+1)
\sum\limits_{k_1=0}^
{h_{(s-1)}-1}\sum\limits_{k_2=0}^
{h_{(s-2)}-h_{(s-1)}}\ldots\sum\limits_{k_{(s-2)}=0}^{h_2-h_3}
(r-h_1+h_2-2\sum\limits_{i=1}^{s-2}k_i),$}
\item {\small $g(r,h_1,\ldots,h_{(s-1)},1)=
g(r,h_1-h_{(s-1)}+1,\ldots,h_{(s-2)}-h_{(s-1)}+1,1)
+g(r-2,h_1-1,\ldots,h_{(s-1)}-1,1)+
\Delta'(r-2,h_1-h_{(s-1)}+1,
\ldots,h_{(s-2)}-h_{(s-1)}+1,1)=\sum\limits_{k=0}^{h_{(s-1)}-1}
g(r-2k,h_1-h_{(s-1)}+1,\ldots,h_{(s-2)}-h_{(s-1)}+1,1)+
\sum\limits_{k=0}^{h_{(s-1)}-2}
\Delta'(r-2(k+1),h_1-h_{(s-1)}+1,
\ldots,h_{(s-2)}-h_{(s-1)}+1,1).$}
\end{itemize}

Finally, we take $h_s\geq 2$. In
$\b(r,h_1,\ldots,h_s)$, we set that
$\p^{r+h_1}\vee\p^{2r-1}=\p^{2r}$. Then
$\b(r,h_1,\ldots,h_s)$ breaks up into 
$\b(r,h_1-1,\ldots,h_s-1)\sub\p^{2r+1}$ and
$\b(r-2,h_1-1,\ldots,h_s-1)\sub\p^{2r-3}$ with
$\Delta(r-2,h_1-1,\ldots,h_s-1)$ generators in
common. By induction on $h_s\geq1$, we conclude that 
\begin{itemize}
\item {\small
$d(r,h_1,\ldots,h_s)=
\sum\limits_{k=0}^{h_{s}-1}\binom{h_s-1}{k}
d(r-2k,h_1-h_{s}+1,\ldots,h_{(s-1)}-h_{s}+1,1),$}
\item {\small $
g(r,h_1,\ldots,h_s)=
\sum\limits_{k=0}^{h_{s}-1}\binom{h_s-1}{k}
g(r-2k,h_1-h_{s}+1,\ldots,h_{(s-1)}-h_{s}+1,1)+
\sum\limits_{k_1=2}^{h_s}
\sum\limits_{k_2=0}^{h_s-k_1}\binom{h_s-k_1}{k_2}
\Delta'(r-2(k_2+1),
h_1-h_s+k_1-1,\ldots,h_{(s-1)}-h_s+k_1-1,k_1-1).$}
\end{itemize}
\qed

\subsection{Incidence Scrolls with $\e\not\sim
0$}\label{esp}

We can continue with a similar method for provide all the
incidence scrolls which have $\e\not\sim 0$ but $deg(\e)=0$.
Writing $r$ instead of
$m-g+i_1$, there is a one-to-one correspondence between
incidence scrolls in $\p^{2r+1}$ with base $\b=\{2\, \p^{r},
\p^{n_1},\ldots,\p^{n_s}\}$ such that $r<n_1$ and partitions
of $2r-1$ such that $\lambda _1\leq r-1$.
\begin{lemma}\label{delta1} Let $\{h_1,\ldots,h_s\}$ be a 
finite number of positive integers such that
$1\leq h_1\leq r,\, 1\leq h_2\leq r+h_1,\, 1\leq h_i\leq
h_{(i-1)}-1$ for $3\leq i\leq s$ and
$h_s+\sum_{i=3}^{s}(h_{(i-1)}-h_i)=h_2$. Suppose 
$\Delta_{\e\not\sim 0}(r,h_1,\ldots,h_s):=$
$$\left\{\begin{array}{ll}
\displaystyle
\sigma(r)^2\cdot 
\sigma(r+h_1-1)\cdot\sigma(2r-1)^{(r+h_1-1)},&
\mbox{$s=1$}\\
\sigma(r)^2\cdot 
\sigma(r+h_1-1)\cdot\sigma(r-h_1+h_2)\cdot
\sigma(2r-h_2+h_3)\cdot\ldots\cdot& \\
\sigma(2r-h_{(s-1)}+h_s)\cdot
\sigma(2r-1)^{(h_s-1)},& \mbox{$s\geq 2$} 
\end{array} \right.$$ an intersection of Schubert cycles in
$\p^{2r+1}$. Then:
\begin{enumerate}
\item $\Delta_{\e\not\sim
0}(r,h)=\\
=\left\{\begin{array}{ll}
\displaystyle
d_{\e\sim 0}(r,r-1),&h=1\\
\sum\limits_{k_1=h-1}^{r}\sum\limits_{k_2=h-2}^{k_1}
\ldots\sum\limits_{k_{(h-1)}=1}^{k_{(h-2)}}d_{\e\sim
0}(k_{(h-1)},k_{(h-1)}-1),&h>1
\end{array} \right.$
\item $\Delta_{\e\not\sim
0}(r,h_1,h_2)=\sum\limits_{k=0}^{h_2-1}\binom{h_2-1}{k}
(r-h_1+h_2+1-2k);$
\item $\Delta_{\e\sim 0}(r,h_1,\ldots,h_s)=\sum\limits_{k_1=0}^{h_s-1}
\sum\limits_{k_2=0}^{h_{(s-1)}-h_s}
\sum\limits_{k_3=0}^{h_{(s-2)}-h_s-2k_2}\ldots
\sum\limits_{k_{(s-1)}=0}^{h_{2}-h_s-2\sum\limits_{i=2}^{s-2}k_i}
\binom{h_s-1}{k_1}\\(r-h_1+h_s+1-2k_1), \mbox{ for}\;
s\geq3.$
\end{enumerate}
\end{lemma}
{\bf Proof.} We only give the main ideas of the proof
because it is similar to that of Lemma \ref{delta}.
\begin{enumerate}
\item {\small $\Delta(r,1)=d_{\e\sim 0}(r,r-1)$} and 
{\small $\Delta(r,h)=\Delta(r-1,h)+\Delta(r,h-1)=
\sum\limits_{k=h-1}^{r}\Delta(k,h-1)$} (use the fact that
{\small $\Delta(r,r+1)=\Delta(r,r)$}).
\item For $s=2$, we have {\small $\Delta(r,h_1,1)=
\sigma(r-h_1+1)_{_{|2(r-h_1+1)+1}}^4=r-h_1+2$} and
{\small $\Delta(r,h_1,h_2)=\Delta(r,h_1-1,h_2-1)+
\Delta(r-1,h_1,h_2-1)=
\sum\limits_{k_1=0}^{h_2-1}
\binom{h_2-1}{k}\Delta(r-k,h_1-h_2+1+k,1)$}.   
\item When $s\geq3$ and $h_s=1$, we obtain
{\small
$\Delta(r,h_1,\ldots,h_{(s-1)},1)\underset{\mbox{\tiny
Pieri's formula}}{=}\sum\limits_{k=0}^{h_{(s-1)}-1}
\sigma(r-k)^2 
\sigma(r+h_1-1-2k)\sigma(r-h_1+h_2-2k)\sigma(2r-h_2+h_3-2k)\ldots
\sigma(2r-h_{(s-3)}+h_{(s-2)}-2k)_{|_{2(r-k)+1}}=\sum\limits_{k=0}^
{h_{(s-1)}-1}\Delta(r-k,h_1-k,h_2-2k,\ldots,$
$h_{(s-2)}-2k,1)
=\sum\limits_{k_1=0}^{h_{(s-1)}-1}
\sum\limits_{k_2=0}^{h_{(s-2)}-1-2k_1}\ldots
\sum\limits_{k_{(s-2)}=0}^{h_{2}-1-2\sum\limits_{i=1}^{s-3}k_i}
(r-h_1+2).$}

Finally, when $h_s\geq 1$, Pieri's formula makes it
obvious that
{\small $\Delta(r,h_1,\ldots,h_{s})$ $=\sum\limits_{k=0}^
{h_{s}-1}\binom{h_s-1}{k}\Delta
(r-k,h_1-h_s+1+k,h_2-h_s+1,\ldots,h_{(s-1)}-h_s+1,1).$}\qed
\end{enumerate}
\begin{prop}\label{unoseis} Each partition
$(r-h,1^{(r+h-1)})$ of
$2r-1$ such that
$h\geq 1$ gives rise to $R^{d_{\e\not\sim 0}(r,h)}_{g_{\e\not\sim
0}(r,h)}\sub\p^{2r+1}$ with base
$\b_{\e\not\sim
0}(r,h)=\{2\,\p^r,$ $\p^{r+h}, (r+h-1)\, \p^{2r-1}\}$ in
general position, where:
\begin{itemize}
\item $d_{\e\not\sim
0}(r,h)=\sum\limits_{k_1=h}^{r}
\sum\limits_{k_2=h-1}^{k_1}\ldots
\sum\limits_{k_{(h-1)}=2}^{k_{(h-2)}}
\sum\limits_{k_h=1}^{k_{(h-1)}}
d_{\e\sim
0}(k_h,k_h-1)$, 
\item $g_{\e\not\sim 0}(r,h)=\sum\limits_{k_1=h}^{r}
\sum\limits_{k_2=h-1}^{k_1}\ldots
\sum\limits_{k_{(h-1)}=2}^{k_{(h-2)}}\sum\limits_{k_h=1}^{k_{(h-1)}}
g_{\e\sim0}(k_h,k_h-1)+\sum\limits_{\A=2}^{h}\Bigl[
\sum\limits_{k_1=h}^{r}\\\sum\limits_{k_2=h-1}^{k_1}\ldots
\sum\limits_{k_{(\A-1)}=h-\A+2}^{k_{(\A-2)}}
\sum\limits_{k_{\A}=h-\A+1}^{k_{(\A-1)}-1}\Delta'_{\e\not\sim
0}(k_{\A},h-\A+1)\Bigr]+\sum\limits_{k=h}^{r-1}\Delta'_{\e\not\sim
0}(k,h)$.
\end{itemize}
\end{prop}
{\bf Proof}. We proceed by induction on $h\geq1$. For
$h=1$, suppose that, in $\b(r,1)$,
$\p^{r+1}\vee\p^{2r-1}=\p^{2r}$. Then $\b(r,1)$ degenerates
into $\b_{\e\sim 0}(r,r-1)\sub\p^{2r+1}$ and
$\b_{\e\not\sim 0}(r-1,1)\sub \p^{2r-1}$
with $\Delta_{\e\not\sim 0}(r-1,1)$ generators in common.
Whence, by induction on $r\geq 2$, we find:
\begin{itemize}
\item $d(2,1)=d(2,2)=d_{\e\sim
0}(2,1)+d_{\e\sim 0}(1,0)$ and
$d(r,1)=d_{\e\sim 0}(r,r-1)+d(r-1,1)=
\sum\limits_{k=1}^{r}d_{\e\sim 0}(k,k-1),$
\item $g(2,1)=g_{\e\sim
0}(2,1)+g_{\e\sim 0}(1,0)+\Delta'(1,1)$ and
$g(r,1)=g_{\e\sim 0}(r,r-1)+g(r-1,1)+\Delta'(r-1,1)=
\sum\limits_{k=1}^{r}g_{\e\sim
0}(k,k-1)+\sum\limits_{k=1}^{r-1}\Delta'(k,1)
.$
\end{itemize}

Assume the formulas hold for $h-1$; we will prove them for
$h>1$. If $\p^{r+h}\vee\p^{2r-1}=\p^{2r}$, then $\b(r,h)$
breaks up into 
$\b(r,h-1)\sub\p^{2r+1}$ and $\b(r-1,h)\sub\p^{2r-1}$ with
$\Delta_{\e\not\sim 0}(r-1,h)$ generators in common.
Therefore,
\begin{itemize}
\item $d(r,h)=d(r,h-1)+d(r-1,h)=\sum\limits_{k=h}^{r}d(k,h-1),$
\item $g(r,h)=g(r,h-1)+g(r-1,h)+\Delta'(r-1,h)
=\sum\limits_{k=h}^{r}g(k,h-1)
+\sum\limits_{k=h}^{r-1}\Delta'(k,h).$
\end{itemize}
By induction on $r\geq h$, since $\b(r,r)=\b(r,r-1)$, the
last equalities are obvious. Combining the two results,
we can rewrite:
\begin{itemize}
\item $d(r,h)=
\sum\limits_{k_1=h}^{r}
\sum\limits_{k_2=h-1}^{k_1}\ldots
\sum\limits_{k_{(h-1)}=2}^{k_{(h-2)}}d(k_{(h-1)},1)$,
\item $g(r,h)=
\sum\limits_{k_1=h}^{r}
\sum\limits_{k_2=h-1}^{k_1}
\ldots
\sum\limits_{k_{(h-1)}=2}^{k_{(h-2)}}
g(k_{(h-1)},1)+\sum\limits_{k=h}^{r-1}\Delta'(k,h)+\sum\limits_{\A=2}^{h-1}\Bigl[
\sum\limits_{k_1=h}^{r}\\
\sum\limits_{k_2=h-1}^{k_1}\ldots
\sum\limits_{k_{(\A-1)}=h-\A+2}^{k_{(\A-2)}}
\sum\limits_{k_{\A}=h-\A+1}^{k_{(\A-1)}-1}\Delta'(k_{\A},h-\A+1)\Bigr].$
\end{itemize}
\begin{teo} \label{eno} For $s\geq2$, each partition
$(r-h_1,r+h_1-h_2,h_2-h_3\ldots,h_{(s-1)}-h_s,1^{(h_s-1)})$
of
$2r-1$ such that
$h_1\geq 1$
gives rise to the scroll $R^{d_{\e\not\sim
0}(r,h_1,\ldots,h_{s})}_{g_{\e\not\sim
0}(r,h_1,\ldots,h_{s})}$ $\sub\p^{2r+1}$ with base
$\b_{\e\not\sim 0}(r,h_1,\ldots,h_s)=\{2\,\p^r, 
\p^{r+h_1},\p^{r-h_1+h_2},$ $\p^{2r-h_2+h_3},$ $\ldots,
\p^{2r-h_{(s-1)}+h_s}, (h_s-1)\,\p^{2r-1}\}$ in
general position, where
\begin{itemize}
\item {\scriptsize $d_{\e\not\sim 0}(r,h_1,\ldots,h_{s})=
2\sum\limits_{k_1=0}^{h_s-1}\sum\limits_{k_2=0}^
{h_{(s-1)}-h_s}\ldots\sum\limits_{k_{(s-1)}=0}^{h_2-h_3}
\binom{h_s-1}{k_1}
(r-h_1+h_2
-2\sum\limits_{i=1}^{s-1}k_i),$}
\item {\scriptsize $g_{\e\not\sim
0}(r,h_1,\ldots,h_{s})=
\sum\limits_{\A=1}^{s-2}\Bigl[\sum\limits_{k_1=0}^{h_s-1}
\sum\limits_{k_2=0}^
{h_{(s-1)}-h_s}\ldots
\sum\limits_{k_{\A}=0}^{
h_{(s-\A+1)}-h_{(s-\A+2)}}
\sum\limits_{k_{(\A+1)}=0}^
{h_{(s-\A)}-h_{(s-\A+1)}-1}\\
\binom{h_s-1}{k_1}\delta'(r-1-
\sum\limits_{i=1}^{\A+1}k_i,h_1-h_{(s-\A)}
+1+\sum\limits_{i=1}^{\A+1}k_i,h_2-h_{(s-\A)}+1,\ldots,
h_{(s-\A-1)}-h_{(s-\A)}\\+1,1)\Bigr]+
\sum\limits_{k_1=0}^{h_s-2}\sum\limits_{k_2=0}^
{k_1}\binom{k_1}{k_2}\Delta'_{\e\not\sim
0}(r-k_2-1,h_1-k_1+k_2,h_2-1-k_1,\ldots,h_s-1-k_1);
$}
\end{itemize}
and, for $s\geq3$,  
{\footnotesize
$\delta(r,h_1,\ldots,h_{s-2},1)=
\sum\limits_{k_1=0}^{h_{(s-2)}-1}\sum\limits_{k_2=0}^
{h_{(s-3)}-1-2k_1}\ldots
\sum\limits_{k_{(s-3)}=0}^{h_{2}-1-2
\sum\limits_{i=1}^{s-4}k_i}(r-h_1+1)=
\delta'(r,h_1,\ldots,h_{s-2},1)+1$}.
\end{teo}
{\bf Proof.} First, we will compute the case $s=2$. For
$h_2=1$, we obtain $\b(r,h_1,1)=\{3\,\p^{r-h_1+1},
\p^{r-h_1+2}\}$, hence $d(r,h_1,1)=2(r+h_1+1)$ and
$g(r,h_1,1)=0$. For
$s\geq3$ and
$h_s=1$, we can suppose that
$\p^{r+h_1}\cap\p^{2r-h_{(s-1)}+1}$ $=\p^{r+h_1-h_{(s-1)}+1}$.
Then the scroll degenerates into
\begin{enumerate}
\item $\b(r,h_1-h_{(s-1)}+1,\ldots,h_{(s-2)}-h_{(s-1)}+1,1)
\sub\p^{2r+1}$,
\item
$\b(r-1,h_1,h_2-1,\ldots,h_{(s-1)}-1,1)
\sub\p^{2r-1}$, with
$\sigma(r-1)^2\sigma(r+h_1-h_{(s-1)})\sigma(r-h_1+h_2-2)
\sigma(2r-h_2+h_3-2)\ldots\sigma(2r-h_{(s-2)}+h_{(s-1)}-2)_{|_{2r-1}}$
generators in common.
\end{enumerate}
We set
$\delta(r,h_1,\ldots,h_{(s-2)},1):=\sigma(r)^2\sigma(r+h_1)
\sigma(r-h_1+h_2-1)
\sigma(2r-h_2+h_3)\ldots\sigma(2r-h_{(s-3)}+h_{(s-2)})
\sigma(2r-h_{(s-2)}+1)_{|_{2r+1}}$. From
Pieri's formula,
$$\delta(r,h_1,\ldots,h_{(s-2)},1)=\sum\limits_{k=0}^{h_{(s-2)}-1}
\delta(r-k,h_1-k,h_2-2k,\ldots,h_{(s-3)}-2k,1)$$ 
Our assertion follows by induction on $s$ because it is clear
that  
$\delta(r,h_1,1)=r-h_1+1$. Then \begin{itemize}
\item {\small
$d(r,h_1,\ldots,h_{(s-1)},1)=\sum\limits_{k=0}^{h_{(s-1)}-1}
d(r-k,h_1-h_{(s-1)}+1+k,h_2-h_{(s-1)}+1,\ldots,h_{(s-2)}-h_{(s-1)}+1,
1)$,}
\item {\small
$g(r,h_1,\ldots,h_{(s-1)},1)=\sum\limits_{k=0}^{h_{(s-1)}-1}
g(r-k,h_1-h_{(s-1)}+1+k,h_2-h_{(s-1)}+1,\ldots,h_{(s-2)}
-h_{(s-1)}+1, 1)+\sum\limits_{k=0}^{h_{(s-1)}-2}
\delta'(r-k-1,h_1-h_{(s-1)}+1+k,h_2-h_{(s-1)}+1,\ldots,h_{(s-2)}-h_{(s-1)}
+1,1).$}
\end{itemize}
By induction on $s\geq2$, we obtain:
\begin{itemize}
\item{\small
$d(r,h_1,\ldots,h_{(s-1)},1)=2\sum\limits_{k_1=0}^{h_{(s-1)}-1}
\sum\limits_{k_2=0}^
{h_{(s-2)}-h_{(s-1)}}\ldots
\sum\limits_{k_{(s-2)}=0}^{h_{2}-h_3}(r-h_1+h_2-2\sum\limits_{i=1}^{s-2}
k_i),$}
\item{\small
$g(r,h_1,\ldots,h_{(s-1)},1)=\sum\limits_{k=0}^{h_{(s-1)}-2}
\delta'(r-k-1,h_1-h_{(s-1)}+1+k,h_2-h_{(s-1)}+1,\ldots,h_{(s-2)}-h_{(s-1)}
+1,1)+\sum\limits_{\A=2}^{s-2}\Bigl[\sum\limits_{k_1=0}^{h_{(s-1)}-1}
\sum\limits_{k_2=0}^
{h_{(s-2)}-h_{(s-1)}}\ldots
\sum\limits_{k_{(\A-1)}=0}^{h_{(s-\A+1)}-h_{(s-\A+2)}}\\
\sum\limits_{k_{\A}=0}^{h_{(s-\A)}-h_{(s-\A+1)}-1}
\delta'(r-1-\sum\limits_{i=1}^{\A}
k_i,h_1-h_{(s-\A)}+1+\sum\limits_{i=1}^{\A}
k_i,h_2-h_{(s-\A)}+1,\ldots,\\h_{(s-\A-1)}-h_{(s-\A)}
+1,1)\Bigr].$}
\end{itemize}

Finally, take $h_s\geq2$. In $\b(r,h_1,\ldots,h_s)$, suppose that
$\p^{r+h_1}\cap\p^{2r-1}=\p^{r+h_1-1}$. Then the scroll
degenerates into
$\b(r,h_1-1,\ldots,h_s-1)\sub\p^{2r+1}$ and
$\b(r-1,h_1,h_2-1,\ldots,h_s-1)\sub\p^{2r-1}$ with
$\Delta_{\e\not\sim
0}(r-1,h_1,h_2-1,\ldots,h_s-1)$ generators in
common. For this, 
\begin{itemize} 
\item{\small $d(r,h_1,\ldots,h_s)=
\sum\limits_{k=0}^{h_{s}-1}\binom{h_{s}-1}{k}
d(r-k,h_1-h_{s}+1+k,h_2-h_{s}+1,\ldots,h_{(s-1)}-h_{s}+1,1)$,}
\item{\small $g(r,h_1,\ldots,h_s)=
\sum\limits_{k=0}^{h_{s}-1}\binom{h_{s}-1}{k}
g(r-k,h_1-h_{s}+1+k,h_2-h_{s}+1,\ldots,h_{(s-1)}-h_{s}+1,1)+
\sum\limits_{k_1=0}^{h_{s}-2}
\sum\limits_{k_2=0}^{k_1}\binom{k_1}{k_2}
\Delta'_{\e\not\sim
0}(r-k_2-1,h_1-k_1+k_2,h_2-1-k_1,\ldots,h_{s}-1-k_1)$.}\qed
\end{itemize}

In case $s=2$, if $(r-h_1,r+h_1-h_2,1^{(h_2-1)})$ is a
partition of
$2r-1$, then it is required that $2h_1\leq h_2$. But the
above formulas are true for all positive integers $h_1,
h_2$. This is possible because $\b(r,h_1,h_2)$ is base of an
incidence scroll although $(r-h_1,r+h_1-h_2,1^{(h_2-1)})$ is
not a partition.

In the general case , all incidence scroll
with $\e\not\sim0$ is one $R^d_g\sub \p^{2r-e+1+j}$ with
$\{\p^{r-e}, \p^{r+j}\}\sub\b$ ($r:=m-g+i_1$ and
$j:=i_2-i_1$). Moreover, there are a one-to-one correspondence
between partitions of $2r+e+j+1$ such that $\lambda _1\leq
r-e$ and bases of such scrolls.

\begin{lemma}\label{delta2} Let $\{h_1,\ldots,h_s\}$ be a
finite number of positive integers such that 
$0\leq h_1\leq r-e-1,\, 1\leq h_2\leq r+j+h_1,\,1\leq
h_i\leq h_{(i-1)}-1$ for $3\leq i\leq s$ and
$h_s+\sum_{i=3}^{s}(h_{(i-1)}-h_i)=h_2$. Define
$\Delta_{e\geq 1}(r,h_1,\ldots,h_s):=$
$$\left\{\begin{array}{ll}
\displaystyle
\sigma(r-e)\cdot\sigma(r+j)\cdot 
\sigma(r+j+h_1-1)\cdot\sigma(2r-e+j-1)^{(r+j+h_1-1)},& \mbox{$s=1$}\\
\sigma(r-e)\cdot\sigma(r+j)\cdot
\sigma(r+j+h_1-1)\cdot\sigma(r-e-h_1+h_2)\cdot&\\
\sigma(2r-e+j-h_2+h_3)\cdot\ldots\cdot
\sigma(2r-e+j-h_{(s-1)}+h_s)\cdot& \\
\sigma(2r-e+j-1)^{(h_s-1)},& \mbox{$s\geq 2$} 
\end{array} \right.$$ an intersection of
Schubert cycles in $\p^{2r-e+j+1}$. Then:
\begin{enumerate}
\item
$\Delta_{e\geq
1}(r,h)=\sum\limits_{k_1=0}^{r-e-h+1}\sum\limits_{k_2=0}^{r-e-h+2-k_1}
\ldots\sum\limits_{k_{h}=0}^{r-e-\sum\limits_{j=1}^{h-1}k_j}d_{e
\geq 1}(r-\sum\limits_{j=1}^{h}k_j,0)_{(j-1)};$
\item $\Delta_{e\geq
1}(r,h_1,1)=\left \{\begin{array}{ll}
\displaystyle
r-e+1 & \mbox{if $h_1=0$}\\
r-e-h_1+2 & \mbox{if $h_1>0$}
\end{array} \right. ;$
\item $\Delta_{e\geq
1}(r,h_1,\ldots,h_s)=\sum\limits_{k_1=0}^{h_s-1}
\sum\limits_{k_2=0}^{h_{(s-1)}-h_s}
\sum\limits_{k_3=0}^{h_{(s-2)}-h_s-2k_2}\ldots
\sum\limits_{k_{(s-1)}=0}^{h_{2}-h_s-2\sum\limits_{j=2}^{s-2}k_j}
\binom{h_s-1}{k_1}\\
\Delta_{e\geq
1}(r-\sum\limits_{j=2}^{s-1}k_j,h_1-h_s+1+k_1
-\sum\limits_{j=2}^{s-1}k_j,1),
 \; s\geq3.$
\end{enumerate}
 \end{lemma}

The number $d_{e\geq
1}(r,0)_{j}$ will be computed in Proposition
\ref{unonueve}.  

{\bf Proof.} If $s=1$, then {\small 
$\Delta(r,0)=\sigma(r-e)\sigma(r+j-1)^2\sigma
(2r-e+j-2)^{r+j-1}_{|_{2r-e+j}}=d_{e\geq1}(r,0)_{j-1}$} (see
Proposition \ref{unonueve}). In general ($h\geq 1$),
apply Pieri's formula to $\sigma(r+j+h-1)$ and any
$\sigma(2r-e+j-1)$. Then, by induction on $r\geq e+h+1$,
$$\Delta(r,h)=\Delta(r-1,h)+\Delta(r,h-1)=
\sum\limits_{k=0}^{r-e-h+1}\Delta(r-k,h-1).$$ Use that 
$\Delta(e+h-1,h)=\Delta(e+h-1,h-1)$.

For $s=2$, $\Delta(r,0,1)=r-e+1$
because is the degree of the directrix curve, which is contained
in an $\p^{r-e+1}$, of $R^{2(r-e)+1}_0\p^{2(r-e)+2}$; and
$\Delta(r,h_1\geq1,1)=\sigma(r-e-h_1+1)^4_{|_{2(r-e-h_1+1)+1}}
=r-e-h_1+2$.

The other cases are left to the reader
because these are analogous of those of Lemma \ref{delta}.
\qed
\begin{prop}\label{unonueve} Each partition
$(r-e-h,1^{(r+j+h-1)})$ of
$2r-e+j-1$ such that $h\geq 0$ give rise to
$R^{d_{e\geq 1}(r,h)}_{g_{e\geq 1}(r,h)}\sub
\p^{2r-e+j+1}$ with base $\b_{e\geq1}(r,h)=\{\p^{r-e},
\p^{r+j},
\p^{r+j+h}, (r+j+h-1)\, \p^{2r-e+j-1}\}$ in
general position, where
\begin{itemize}
\item {\scriptsize $d_{e\geq
1}(r,h)=
\sum\limits_{k_1=0}^{r-e-h}
\sum\limits_{k_2=0}^{r-e-h+1-k}\ldots
\sum\limits_{k_{(h-1)}=0}^{r-e-2-\sum\limits_{i=1}^{h-2}k_i}
\sum\limits_{k_h=0}^{r-e-1-\sum\limits_{i=1}^{h-1}k_i}\Bigl[1+
\sum\limits_{k_{(h+1)}=0}^{r-e-1-\sum\limits_{i=1}^{h}k_i}\\
\sum\limits_{k_{(h+2)}=0}^{r-e-1-
\sum\limits_{i=1}^{h+1}k_i}
\ldots
\sum\limits_{k_{(h+e+j-1)}=0}^{r-e-1-\sum\limits_{i=1}^{h+e+j-2}k_i}
d_{\e\not\sim0}(r-e-\sum\limits_{i=1}^{h+e+j-1}k_i,1)+
\sum\limits_{\A=0}^{j+e-2}
\bigl[\sum\limits_{k_{(h+1)}=0}^{r-e-1-\sum\limits_{i=1}^{h}k_i}\ldots\\
\sum\limits_{k_{(h+\A)}=0}^{r-e-1-\sum\limits_{i=1}^{h+\A-1}k_i}(r-e-
\sum\limits_{i=1}^{h+\A}k_i)\bigr]\Bigr]$,}
\item{\small $g_{e\geq
1}(r,h)=\sum\limits_{k_1=0}^{r-e-h}
\sum\limits_{k_2=0}^{r-e-h+1-k}\ldots
\sum\limits_{k_{(h-1)}=0}^{r-e-2-\sum\limits_{i=1}^{h-2}k_i}
\sum\limits_{k_h=0}^{r-e-1-\sum\limits_{i=1}^{h-1}k_i}\Bigl[
\sum\limits_{k_{(h+1)}=0}^{r-e-2-\sum\limits_{i=1}^{h}k_i}
\ldots\\
\sum\limits_{k_{(h+e+j-1)}=0}^{r-e-2-\sum\limits_{i=1}^{h+e+j-2}k_i}
g_{\e\not\sim0}(r-e-\sum\limits_{i=1}^{h+e+j-1}k_i,1)+
\sum\limits_{\A=1}^{j+e-1}\bigl[
\sum\limits_{k_{(h+1)}=0}^{r-e-2-\sum\limits_{i=1}^{h}k_i}\ldots\\
\sum\limits_{k_{(h+\A)}=0}^{r-e-2-\sum\limits_{i=1}^{h+\A-1}k_i}\beta
(r-1-
\sum\limits_{i=1}^{h+\A}k_i,j-\A)\bigr]\Bigr]+
\sum\limits_{k=0}^{r-e-h-1}\Delta'_{e\geq
1}(r-k-1,h)+
\sum\limits_{\A=2}^{h-1}\Bigl[
\sum\limits_{k_1=0}^{r-e-h}\\
\sum\limits_{k_2=0}^{r-e-h+1-k}\ldots
\sum\limits_{k_{(\A-1)}=0}^{r-e-h+\A-2-\sum\limits_{i=1}^{\A-2}k_i}
\sum\limits_{k_{\A}=0}^{\;
r-e-h+\A-2-\sum\limits_{i=1}^{\A-1}k_i}
\Delta'_{e\geq
1}(r-1-\sum\limits_{i=1}^{\A}k_i,h-\A+1)\Bigr]$,}
\end{itemize}
with
$\beta(r,j)=d(r,0)_j-1$.
\end{prop}
{\bf Proof}. We proceed by induction on $h\geq0$. 
For $h=0$, let us the temporary notation of $(r,0)_j$ for
$(r,0)$. It is clear that $r+h+j>1$.
In  other case, we have $\b=\{3\p^1\}\sub \p^3$ but
$e\not\sim0$. In $\b(r,0)_j$, suppose that
$\p^{r-e}\cap\p^{2r-e+j-1}=\p^{r-e-1}$. Then the scroll
breaks up into $\b(r,0)_{j-1}\sub\p^{2r-e+j}$ and
$\b(r-1,0)_j\sub
\p^{2r-e+j-1}$ with 
$\sigma(r-e-1)\sigma(r+j-1)^2\sigma(2r-e+j-2)^{r+j-2}
=\sigma(r-e-1)\sigma(r+j-2)^2\sigma(2r-e+j-4)^{r+j-2}_{|_{2r-e+j-2}}
=d(r-1,0)_{j-1}$ generators in common. By induction on $r\geq
e+1$, we find that
\begin{itemize}
\item $d(r,0)_j=1+\sum\limits_{k=0}^{r-e-1}d(r-k,0)_{j-1},$
\item $g(r,0)_j=\sum\limits_{k=0}^{r-e-2}g(r-k,0)_{j-1}+
\sum\limits_{k=0}^{r-e-2}[d(r-1-k,0)_{j-1}-1].$
\end{itemize}
When $j=-e+1$, it is easy to check that 
$d(r,0)_{(-e+1)}=1+d_{\e\not\sim 0}(r-e,1)$ and
$g(r,0)_{(-e+1)}=g_{\e\not\sim 0}(r-e,1)$. Hence,
by induction on $j\geq -e+1$, 
\begin{itemize}
\item $d(r,0)_j=1+\sum\limits_{k_1=0}^{r-e-1}
\sum\limits_{k_2=0}^{r-e-1-k}\ldots
\sum\limits_{k_{(j+e-1)}=0}^
{r-e-1-\sum\limits_{i=1}^{j+e-2}k_i}d_{\e\not\sim
0}(r-e-\sum\limits_{i=1}^{j+e-1}k_i,1)+
\sum\limits_{\A=0}^{j+e-2}\Bigl[
\sum\limits_{k_1=0}^{r-e-1}
\sum\limits_{k_2=0}^{r-e-1-k}\ldots
\sum\limits_{k_{\A}=0}^
{r-e-1-\sum\limits_{i=1}^{\A-1}k_i}(r-e-\sum\limits_{i=1}^{\A}k_i)
\Bigr]$,
\item $g(r,0)_j=\sum\limits_{k_1=0}^{r-e-2}
\sum\limits_{k_2=0}^{r-e-2-k}\ldots
\sum\limits_{k_{(j+e-1)}=0}^
{r-e-2-\sum\limits_{i=1}^{j+e-2}k_i}g_{\e\not\sim
0}(r-e-\sum\limits_{i=1}^{j+e-1}k_i,1)+
\sum\limits_{\A=1}^{j+e-1}\\
\Bigl[
\sum\limits_{k_1=0}^{r-e-2}
\sum\limits_{k_2=0}^{r-e-2-k}\ldots
\sum\limits_{k_{\A}=0}^
{r-e-2-\sum\limits_{i=1}^{\A-1}k_i}
\beta(r-1-\sum\limits_{i=1}^{\A}k_i,j-\A)
\Bigr]$.
\end{itemize}
Assuming the formulas hold for $h-1\geq0$, we will prove
them for $h$. In $\b(r,h)$, suppose that
$\p^{r+j+h}\cap\p^{2r-e+j-1}=\p^{r+j+h-1}$. Then $R$
degenerates into
$\b(r,h-1)\sub\p^{2r-e+j+1}$ and $\b(r-1,h)\sub\p^{2r-e+j-1}$
with $\Delta_{e\geq1}(r-1,h)$ generators in common. 
By induction on $r\geq e+h+1$, we find that 
\begin{itemize}
\item{\small $d(r,h)=\sum\limits_{k=0}^{r-e-h}d(r-k,h-1)=
\sum\limits_{k_1=0}^{r-e-h}
\sum\limits_{k_2=0}^{r-e-h+1-k}\ldots
\sum\limits_{k_{h}=0}^{r-e-1-\sum\limits_{i=1}^{h-1}k_i}
d(r-\sum\limits_{i=1}^{h}k_i,0)$,}
\item{\small $g(r,h)=\sum\limits_{k=0}^{r-e-h}g(r-k,h-1)+
\sum\limits_{k=0}^{r-e-h-1}\Delta'(r-k-1,h)=\sum\limits_{k_1=0}^{r-e-h}
\sum\limits_{k_2=0}^{r-e-h+1-k}\ldots\\
\sum\limits_{k_{h}=0}^{r-e-1-\sum\limits_{i=1}^{h-1}k_i}
g(r-\sum\limits_{i=1}^{h}k_i,0)+
\sum\limits_{\A=2}^{h-1}\Bigl[
\sum\limits_{k_1=0}^{r-e-h}
\sum\limits_{k_2=0}^{r-e-h+1-k}\ldots
\sum\limits_{k_{(\A-1)}=0}^{r-e-h+\A-2-\sum\limits_{i=1}^{\A-2}k_i}
\\
\sum\limits_{k_{\A}=0}^{\;
r-e-h+\A-2-\sum\limits_{i=1}^{\A-1}k_i}
\Delta'(r-1-\sum\limits_{i=1}^{\A}k_i,h-\A+1)\Bigr]+
\sum\limits_{k=0}^{r-e-h-1}\Delta'(r-k-1,h);$}
\end{itemize}
The last equalities can be proved by induction on $h\geq0$.\qed

\begin{rem}{\em For $\b_{e\geq
1}(r,h_1,h_2)=\{\p^{r-e},\p^{r+j}, 
\p^{r+j+h_1},\p^{r-e-h_1+h_2},
(h_2-1)\p^{2r-e+j-1}\}$ in general position, it follows
easily that
$$d_{e\geq 1}(r,h_1,1)=\left
\{\begin{array}{ll}
\displaystyle
2(r-e)+1 & \mbox{if $h_1=0$}\\
2(r-e-h_1)+2 & \mbox{if $h_1>0$}
\end{array} \right. ;\quad
g_{e\geq 1}(r,h_1,1)=0.$$ 
If $h_2=1$, then
$e+j+2h_1\leq 1$ but $e+j>0$, i.e.,
$h_1=0$. For this $d(r,0,1)=2(r-e)+1$ and $g(r,0,1)=0$. In
case $h_1\geq 1$,we obtain
$\b(r,h_1,1)=\{3\,\p^{r-e-h_1+1},\p^{r-e-h_1+2}\}
\sub\p^{2(r-e-h_1+1)+1}$, i.e., $d(r,h_1,1)=2(r-e-h_1+1)$ and
$g(r,h_1,1)=0$. This result will be needed in the following
theorem which may be proved in much the same way as Theorem
\ref{eno}}
\end{rem}

\begin{teo}\label{em1} For $s\geq2$, each partition
$(r-e-h_1,r+j+h_1-h_2,h_2-h_3,\ldots,h_{(s-1)}-h_s,1^{(h_s-1)})$
of $2r-e+j-1$ such that
$h_1\geq 0$
corresponds to the scroll $R^{d_{e\geq
1}(r,h_1,\ldots,h_{s})}_{g_{e\geq
1}(r,h_1,\ldots,h_{s})}\sub\p^{2r-e+j+1}$ with base  
$\b_{e\geq 1}(r,h_1,\ldots,h_{s})=\{\p^{r-e},$ $\p^{r+j}, 
\p^{r+j+h_1},\p^{r-e-h_1+h_2},\p^{2r-e+j-h_2+h_3},\ldots,\p^{2r-e+j-h_{(s-1)}+h_s},
(h_s-1)\quad$ $\p^{2r-e+j-1}\}$ in general position, where
\begin{itemize}
\item{\small $d_{e\geq 1}(r,h_1,\ldots,h_{s})=
\sum\limits_{k_1=0}^{h_s-1}\sum\limits_{k_2=0}^
{h_{(s-1)}-h_s}\ldots\sum\limits_{k_{(s-1)}=0}^{h_2-h_3}\binom{h_s-1}{k_1}
d(r-\sum\limits_{i=1}^{s-1}
k_i,h_1-h_2+1+\sum\limits_{i=1}^{s-1}k_i,1),$}
\item{\small $g_{e\geq 1}(r,h_1,\ldots,h_{s})=
\sum\limits_{\A=1}^{s-2}\Bigl[\sum\limits_{k_1=0}^{h_s-1}
\sum\limits_{k_2=0}^
{h_{(s-1)}-h_s}\ldots\sum\limits_{k_{\A}=0}^{h_{(s-\A+1)}-h_{(s-\A+2)}}
\sum\limits_{k_{(\A+1)}=0}^{h_{(s-\A)}-h_{(s-\A+1)}-1}\\
\binom{h_s-1}{k_1}\bigl[\delta(r-1-\sum\limits_{i=1}^{\A+1}k_i,h_1-
h_{(s-\A)}+1+\sum\limits_{i=1}^{\A+1}k_i,
h_2-
h_{(s-\A)}+1,h_{(s-\A-1)}-h_{(s-\A)}+1,1)-1\bigr]\Bigr]+
\sum\limits_{k_1=0}^{h_s-2}\sum\limits_{k_2=0}^
{k_1}\binom{k_1}{k_2}
\Delta'_{e\geq1}
(r-k_2-1,
h_1-k_1+k_2,h_2-1-k_1,\ldots,h_{s}-1-k_1),$}
\end{itemize}
letting
$\beta'(r,h_1,\ldots,h_{(s-2)},1)=\delta'(r-e,h_1,\ldots,h_{(s-2)},1).$\qed
\end{teo}

\begin{rem}{\em In general, 
$\b_{e\geq1}(r,h_1,\ldots,h_{s})\sub \p^{2r-e+j+1}$ is not
 base of an incidence scroll which corresponds to an
abstract model with $deg(-\e)$ $=e$. For example, for $s=1,
\; h_1=0$ and
$r=e+1$, we know that $\b(e+1,0)$ generates
$R^{e+j+2}_0\sub\p^{e+j+3}$. It is the 
immersion of the rational ruled surface
$X_{(e+j)}=\p(\Te_{\p^1}\oplus\Te_{\p^1}(-e-j))$ by the very
ample divisor $H\sim C_o+(e+j)f$, i.e., $deg(-\e)=e+j\geq1$.
For this reason, if we talk about $\p^{2r-e+j+1}$, then we
need only consider the numbers
$r,e$ and $j$ such that
$r+j>r-e\geq1$. So we obtain a decomposable incidence scroll
with $deg(-\e)\geq1$.}
\end{rem}

\begin{teo} Let $\pi:X=\p(\E)\lrw C$ be a decomposable
ruled surface over a curve $C$ of genus $g$, such that
$h^0(\Te_X(C_0-\e f))=1$. Let
$\tilde X_x=\p(\tilde {\E})$ be the elementary transform of
$X$ at $x$ where $x\not\in C_o$ and $x\not\in C_1$. Then, 
every projective model of $X$ as an incidence scroll
$R^d_g\sub \p^n$ with base
$\b=\{\p^{n_1},\p^{n_2},\ldots,\p^{n_r}\}$ in general
position, gives rise to an projective model of
$\tilde X_x$, which is an incidence scroll $R^{d-1}_g\sub
\p^{n-1}$ with base
$\tilde{\b}=\{\p^{n_1},\p^{n_2},\p^{n_3-1},\ldots,\p^{n_r-1}\}$
in general position.

Furthermore, these are all the indecomposable incidence
scrolls of $\p^{n-1}$ with $n_1+n_2=n-1$.
\end{teo}
{\bf Proof.} Let $\Phi_{|H|} \colon X
\lrw R^d_g\sub\p^n$ be the immersion induced by $H\sim
C_o+\d f$ on $X$. Let $P=\Phi_{|H|}(x)$ with $y=\pi(x)$.
Then the elementary transform of $X$ at
$x$ corresponds to the projection of $R^d_g\sub \p^n$ from $P$.
Hence
$$\setlength{\unitlength}{5mm}
\begin{picture}(14,3.6)

\put(2.6,3){\makebox(0,0){$X=\p(\E)$}}
\put(2.6,0){\makebox(0,0){$\tilde X_x=\p(\tilde{\E})$}}
\put(10,3){\makebox(0,0){$R^d_g\sub\p^n$}}
\put(10,0){\makebox(0,0){$R^{d-1}_g\sub\p^{n-1}$}}
\put(10,2.3){\vector(0,-1){2}}
\put(2.2,2.3){\vector(0,-1){2}}
\put(4.3,3){\vector(1,0){4}}
\put(4.4,0){\vector(1,0){3.5}}
\put(6.3,4){\makebox(0,0){$\phi_{|C_o+\d f|}$}}
\put(6.3,-0.8){\makebox(0,0){$\phi_{|\xi^{\star}(C_o)+(\d
-y) f|}$}}
\put(1.5,1.5){\makebox(0,0){$\xi_x^{-1}$}}
\put(10.9,1.5){\makebox(0,0){$\pi_{P}$}}

\end{picture}
$$
\smallskip

Moreover, project $R^d_g\sub \p^n$ from $P$, such that $P\notin
C^{m-e}_g\sub\p^{n_1}$ and $P\notin C^{m}_g\sub\p^{n_2}$
for any $C^{m}_g$, is equivalent to join $\p^{n_1}$ and
$\p^{n_2}$. Then, we obtain an incidence scroll $R^{d-1}_g\sub
\p^{n-1}$ with base
$\tilde{\b}$ which is indecomposable because $n_2<n_3$, or
equivalently $h^0(\Te_X(C_0-\e f))=1$.\qed

In this way, we can now compute the degree and genus of
these particular incidence scroll, using the above results
for decomposable incidence scrolls and Proposition
\ref{degen}.

\end{document}